\documentclass[12pt]{article}

\usepackage[T1]{fontenc}
\usepackage[utf8]{inputenc}
\usepackage{authblk}

\usepackage{amssymb,amsthm,amsmath}
\usepackage{amsthm}

\usepackage{setspace}
\usepackage{cases}
\usepackage{dsfont}
\def\dto{\stackrel{{d}}{\longrightarrow}}
\def\ep{\varepsilon}
\def\ds{\displaystyle}
\def\sumn{\sum^n_{i=1}}

\newtheorem{thm}{Theorem}[section]

\newtheorem{theorem}{Theorem}[section]
\newtheorem{corollary}{Corollary}[section]
\newtheorem{lemma}{Lemma}[section]
\newtheorem{remark}{Remark}[section]

\date{}
\title{Another  look at Bootstrapping the Student $t$-statistic}
\author[1]{\small{ \textbf{Mikl\'{o}s Cs\"{o}rg\H{o}}}\thanks{mcsorgo@math.carleton.ca} }
\author[2]{\small{ \textbf{Yuliya V. Martsynyuk }}\thanks{Yuliya.Martsynyuk@ad.umanitoba.ca   } }
\author[1]{\textbf{Masoud M. Nasari}\thanks{mmnasari@math.carleton.ca}    }
\affil[1]{\small{School of Mathematics and Statistics,  Carleton University, Ottawa, ON, Canada}}
\affil[2]{Department of Statistics, University of Manitoba, Winnipeg, MB, Canada   }
\begin{document}
\maketitle
\vspace{-1 cm}
\begin{center}
\emph{Dedicated to the memory of S\'{a}ndor Cs\"{o}rg\H{o}}
\end{center}

\begin{abstract}
{\footnotesize{Let $X, X_1,X_2,\ldots$ be a sequence of i.i.d. random variables with mean $\mu=E X$. Let $\{v_1^{(n)},\ldots,v_n^{(n)}  \}_{n=1}^\infty$ be  vectors of non-negative random variables (weights), independent of the data sequence $\{X_1,\ldots,X_n  \}_{n=1}^\infty$, and put $m_n=\sumn v_i^{(n)}$.  Consider $ X^{*}_1, \ldots, X^{*}_{m_n}$, $m_n\geq 1$, a bootstrap sample, resulting from \emph{re-sampling} or \emph{stochastically re-weighing}  a random sample $X_1,\ldots,X_n$, $n\geq 1$.  Put $\bar{X}_n= \sumn X_i/n$, the original sample mean, and define $\bar{X^*}_{m_n}=\sumn v_i^{(n)} X_i/m_n$, the bootstrap sample mean. Thus, $\bar{X^*}_{m_n}- \bar{X}_n=\sumn ({v_i^{(n)}}/{m_n}-{1}/{n}) X_i$. Put $V_n^{2}=\sumn ({v_i^{(n)}}/{m_n}-{1}/{n})^2$ and let $S_n^{2}$, $S_{m_{n}}^{*^{2}}$ respectively be the the original sample variance and the bootstrap sample variance. The main aim of this exposition is to study the asymptotic behavior of the bootstrapped  $t$-statistics  $T_{m_n}^{*}:= (\bar{X^*}_{m_n}- \bar{X}_n)/(S_n V_n)$ and $T_{m_n}^{**}:= \sqrt{m_n}(\bar{X^*}_{m_n}- \bar{X}_n)/ S_{m_{n}}^{*} $ in terms of \emph{conditioning on the weights} via assuming that, as  $n,m_n\to \infty$, $\max_{1\leq i \leq n}({v_i^{(n)}}/{m_n}-{1}/{n})^2\big/ V_n^{2}=o(1)$ almost surely   or in probability on the probability space of the weights. In  consequence of these  maximum negligibility conditions of the weights, a characterization of the validity of this approach to the bootstrap is obtained as a direct consequence of the Lindeberg-Feller central limit theorem.  This view of justifying the validity of the bootstrap of i.i.d. observables    is believed to be new. The need for it   arises naturally in practice when exploring the nature of information contained in a random sample via re-sampling, for example. Unlike in the theory of weighted bootstrap with exchangeable weights, in this exposition it is not assumed that the components of the vectors of non-negative weights are exchangeable random variables.
 \emph{Conditioning on the data} is also revisited for Efron's bootstrap weights  under conditions on $n,m_n$ as  $n\to \infty $ that differ from requiring $m_n /n$ to be in the interval  $(\lambda_1,\lambda_2)$ with $0< \lambda_1 < \lambda_2 < \infty$ as in Mason and Shao \cite{MasonShao}. Also, the validity of
  the bootstrapped $t$-intervals is established for both approaches to conditioning. Morover, when conditioning on the sample, our results in this regard are new in that they are shown to hold true when $X $ is in the domain of attraction of the normal law (DAN), possibly with infinite variance, while the ones for $E_{X} X^2 <\infty$ when conditioning on the weights are first time results
  \emph{per se}.      }}
\end{abstract}
\noindent
%\begin{keyword}
Keywords: Conditional Central Limit Theorems, Stochastically Weighted Partial Sums, Weighted Bootstrap.
%\end{keyword}
%\begin{keyword}[class=AMS]
%\kwd[Primary ]{62F40, 60K35}
%\kwd{60K35}
%\kwd[; secondary ]{60K35}
%\end{keyword}
% history:
% \received{\smonth{1} \syear{0000}}
%\tableofcontents

\section{Introduction to the approach taken}\label{sect1}

The main objective of the present paper is to address the possibility of   investigating and concluding the validity of bootstrapped partial sums via conditioning on the random weights. The term bootstrapping here will refer to both re-sampling, like Efron's,   and  stochastically re-weighing the data. We show that a direct consequence of the  Lindeberg-Feller  central limit theorem (CLT) as stated  in Lemma \ref{new}, which is also known as the H\'{a}jek-Sid\'{a}k theorem (cf., e.g., Theorem 5.3 in DasGupta \cite{DasGupta}),  is the only required tool to establish the validity of bootstrapped partial sums of independent and identically distributed (i.i.d.) random variables.     As a consequence of  Lemma \ref{new},  Theorem \ref{thm1}  \emph{characterizes}   valid schemes of bootstrap in general, when conditioning on the weights. Accordingly,   the bootstrap weights do not have to be exchangeable in order for the bootstrap scheme to be valid. This is unlike the method of studying the consistency of the generalized   bootstrapped mean that    was initiated by Mason and Newton \cite{Mason and Newton} in terms of conditioning on the sample (cf. their Theorem 2.1 on the thus conditioned asymptotic normality of linear combination of  exchangeable arrays). The latter approach  relies on Theorem 4.1  of  H\'{a}jek \cite{Hajeck} concerning the asymptotic normality of linear rank statistics.
  \par
 We also investigate the validity of Efron's scheme of bootstrap and also that of the scheme of stochastically re-weighing the observations  by  verifying how the respective  bootstrap weights  satisfy the required maximal negligibility  conditions (cf. Corollaries \ref{cor1} and  \ref{cor2}, respectively).

\par
To illustrate  the different nature of the two approaches to conditioning,   we also study Efron's scheme of bootstrap  applied to i.i.d. observations  via conditioning on the data. When doing this, we view a bootstrap partial sum as a randomly weighted sum of centered multinomial random variables. This enables us to derive conditional central limit theorems for these randomly weighted centered  multinomial random variables via  results of Morris \cite{Morris}. The proofs of our Theorems \ref{theorem2} and \ref{theorem3} in this regard  will be seen to be significantly shorter and simpler in comparison to similar results on bootstrapped partial sums when conditioning on the data.

 \par
For throughout use, let $X,X_{1},X_{2},\ldots$ be a sequence of i.i.d. real valued random variables with  mean $\mu:=E(X)$. For a random sample $X_1,\ldots,X_n$, $n\geq 1$,  Efron's  scheme of bootstrap, cf. \cite{Efron},  is a procedure of re-sampling $m_{n}\geq 1$ times with replacement from the original data in such a way that each $X_{i}$, $1\leq i\leq n$, is selected with probability  $\displaystyle{1/n }$ at a time. The resulting sub-sample    will be denoted by   $X_{1}^{*},\ldots,X_{m_{n}}^{*}$, $m_{n}\geq 1$, and is called the bootstrap sample. The  bootstrap partial sum   is a stochastically   re-weighted version of  the original partial sum of $X_1,\ldots,X_n$, i.e.,
\begin{equation}\label{correction1}
\sum_{i=1}^{m_{n}} X^{*}_{i}= \sum_{i=1}^{n} w_{i}^{(n)} X_{i},
\end{equation}
where, $w^{(n)}_{i}:= \#$  \textrm{of \ times\ the\ index}\ $\displaystyle{i}$ \textrm{is\ chosen\ in\ }\ $m_n$\ \textrm{draws \ with \ replacement} \\ $~~~~~~~~~~~~~~~~~~~~ $ \textrm{from} \ $1,\ldots, i,\ldots, n$
 of the indices of $X_{1},\ldots,X_{i},\ldots,X_{n}$.
\begin{remark}\label{remark1}
In view of the preceding definition of $w^{(n)}_{i}$, $1\leq i \leq n$, they form a   row-wise independent triangular array of random variables such that  $\sum_{1\leq i \leq n} w^{(n)}_{i}=m_n$, and for each $n\geq 1$,
$$(w^{(n)}_{1},\ldots,w^{(n)}_{n})\  \substack{d\\=}\ \ multinomial(m_{n};\frac{1}{n},\ldots,\frac{1}{n}),
$$
i.e., a multinomial distribution of size $m_n$ with respective probabilities $1/n$.
Clearly,  for each $n$, $w^{(n)}_{i}$    are independent from the random sample $X_{i}$, $1\leq i \leq n$.
Weights denoted by $w_{i}^{(n)}$ will stand for triangular multinomial random variables in this context throughout.
\end{remark}
The randomly weighted representation of $\sum_{1\leq i \leq m_{n}}X^{*}_{i}$ as in (\ref{correction1}), in turn, enables one to think of bootstrap in a more general way  in which  the scheme of bootstrap is  restricted to neither Efron's  nor to re-sampling in general. In this exposition the term bootstrap will refer   both to \emph{re-sampling}, such as Efron's, as well as to \emph{stochastically re-weighing} the sample. Both of these schemes   of bootstrap can be viewed and treated as \emph{weighted bootstraps}. As such, throughout  this paper, the notation   $v_{i}^{(n)}$, $1\leq i \leq n$, will stand for   \emph{bootstrap weights}  that are to be determined by the scheme of bootstrap in hand. Thus, to begin with, we consider a  sequence $\{ v_{1}^{(n)},\ldots,v_{n}^{(n)} \}_{n\geq 1}$ of vectors of non-negative random weights, independent of the  data sequence $\{X_1,\ldots,X_n \}_{n\geq 1}$, and put $m_n=\sum_{i=1}^n v_{i}^{(n)}$, $m_n \geq 1$. We \emph{do not assume} that  the components of the vectors of the non-negative weights  in hand are exchangeable random variables.

\par
Consider now a  bootstrap sample $X^{*}_{1},\ldots,X^{*}_{m_n}$, $m_n\geq 1$, which is a result of some \emph{weighted bootstrap} via \emph{re-sampling} or \emph{stochastically re-weighing}  the original random sample $X_{1},\ldots,X_{n}$, $n\geq 1$. Define the bootstrap sample mean  $\bar{X}^{*}_{m_n}:=  \sum_{i=1}^{n} v_{i}^{(n)}X_{i}/ m_n$ and the original sample mean $\bar{X}_{n}:=\sum_{i=1}^{n}X_{i}/ n$.  In view of the above setup of bootstrap weights one can readily see that
\begin{eqnarray}
\bar{X}_{m_n}^{*}-\bar{X}_{n}&=& \sum_{i=1}^{n} \big( \frac{v^{(n)}_{i}}{m_{n}}-\frac{1}{n} \big) X_{i} \nonumber\\
&=& \sum_{i=1}^{n} \big( \frac{v^{(n)}_{i}}{m_{n}}-\frac{1}{n} \big) (X_{i}-\mu).\nonumber
\end{eqnarray}
Hence, when studying bootstrapped $t$-statistics via $\{\bar{X}^{*}_{m_n}-\bar{X}_{n} \}_{n\geq 1}  $ in the sequel, it is important to remember that, to begin with,  the latter sequence of statistics has no direct information about the parameter  of interest    $\mu:=E(X)$.

In particular,  in this paper, the  following two general  forms of bootstrapped $t$-statistics   will be considered.
\begin{equation}\label{weightedboot1}
T^*_{m_{n}} =
\frac{ \displaystyle \sum^n_{i=1} \big( \frac{v_{i}^{(n)}}{m_n} - \frac{1}{n} \big) X_i }
{S_{n} \sqrt{\sum^n_{i=1} \big( \frac{v_{i}^{(n)}}{m_n} - \frac{1}{n} \big)^2} },
\end{equation}
\begin{equation}\label{weightedboot2}
T^{**}_{m_{n}} =
\frac{ \displaystyle \sum^n_{i=1} \big( \frac{v_{i}^{(n)}}{m_n} - \frac{1}{n} \big) X_i}
{S^{*}_{m_{n}} / \sqrt{m_{n}} },
\end{equation}
where $S^{2}_{n}$ and $S^{*^2}_{m_n}$ are respectively the original sample variance and the bootstrapped sample variance, i.e.,
$$ S^{2}_{n}=\sum_{1 \leq i \leq n} (X_i-\bar{X}_n)^{2}\big/n $$
and
$$ S^{*^{2}}_{m_n}=\sum_{1 \leq i \leq m_n} ( X^{*}_i-\bar{X}^{*}_{m_n})^{2}\big/m_{n}.$$
\begin{remark}
In this exposition,  both $T^*_{m_{n}}$ and $T^{**}_{m_{n}}$ will be called bootstrapped versions of the well-known Student $t$-statistic
 \begin{equation}\label{t-statistic}
T_{n}:=\frac{\bar{X}_{n}}{S_{n}\big/ \sqrt{n}}=\frac{\sum_{i=1}^n X_i}{S_n \sqrt{n}}.
\end{equation}
\end{remark}
\begin{remark}\label{zeta remark}
In Efron's scheme of bootstrap  $v^{(n)}_{i}=w^{(n)}_{i}$, $1\leq i \leq n$, and  (\ref{weightedboot2}) is seen to be  the well-known Efron bootstrapped $t$-statistic.
 When the parameter of interest is $\mu=E(X)$, Weng \cite{Weng} suggests the use of    $\sum_{i=1}^{n} \zeta_{i}X_{i}/ m_n $, as an estimator of $\mu$, where $\zeta_{i}$ are i.i.d. $Gamma(4,1)$ random variables  which are assumed to be  independent from the random sample $X_i$, $1\leq i \leq n$, and $m_n=\sumn \zeta_i$.  This approach is used in the so-called Bayesian bootstrap (cf., e.g.,  Rubin \cite{Rubin}). This scheme of bootstrap, in a more general form,   shall be viewed   in  Corollary \ref{cor2} below in the context of conditioning on the i.i.d. positive random variables $\{\zeta_1,\ldots,\zeta_n \}$ as specified there.
\end{remark}
\par
The main objective of this exposition is to show that in the presence of the  introduction of the extra randomness, $v_i^{(n)}$, $1\leq i \leq n$,  as a result of re-sampling or re-weighing,  conditional distributions of the  bootstrapped $t$-statistics  $T^{*}_{m_{n}}$ and $T^{**}_{m_{n}}$ will asymptotically coincide with that of    the original $t$-statistic $T_{n}$. In this paper this problem will be studied  by both of the  two  approaches to conditioning in hand.  In Section \ref{sect2}, based on the Lindeberg-Feller CLT   we conclude  a  \emph{characterization} of the asymptotic behavior of the   bootstrapped mean  via  conditioning on the {\em bootstrap weights}, $v_i^{(n)}$, $1\leq i \leq n$, in terms of a manifold  conditional Lindeberg-Feller type CLT     for $T^*_{m_n}$ and $T^{**}_{m_n}$  when $E X^2 < \infty$ (cf. Theorem \ref{thm1}). Then we show  that the validity of Efron's scheme of bootstrap results directly from Theorem \ref{thm1} for both of the latter bootstrapped $t$-statistics when conditioning on $w_{i}^{(n)}$ as in Remark \ref{remark1}  (cf.  Corollary \ref{cor1}). As another example,  in Corollary \ref{cor2}, the  weights $\zeta_{i}/ m_n$, where $\zeta_{i}$ are positive i.i.d.  random variables independent of $\{X_i, 1 \leq i \leq n \}_{n\geq 1}$,   are considered for re-weighing the original sequence. It is shown that  under appropriate moment conditions for $\zeta_i$, the validity of bootstrapping  the $t$-statistic $T_n$ via conditioning on $\zeta_i$, $1\leq i \leq n$, also   follows from Theorem \ref{thm1} for both $T_{m_n}^{*}$ and $T^{**}_{m_n}$ in these terms as well.    In Section \ref{sect3}, we continue the investigation of the limiting conditional distribution of $T^{**}_{m_n}$, but this time  via conditioning on the sample $X_i, 1\leq i \leq n,\ n\geq 1$, and only for Efron's bootstrap scheme, on assuming that $X\in DAN$  (cf. Theorem \ref{theorem2}).
\par
The aim  of weighted bootstrap via conditioning on the bootstrap weights as in Theorem \ref{thm1} is to provide a  scheme of bootstrapping that suites the observations in hand. In other words, it specifies a  method of re-weighing or re-sampling that leads to the same limit as that of  the original $t$-statistic.
This view of justifying the validity of the bootstrap is believed to be new for the two general forms of the bootstrapped Student $t$-statistics $T^{*}_{m_n}$ and $T^{**}_{m_n}$. The need for this approach to the bootstrap in general  arises naturally in practice when exploring the nature of information contained  in a random sample  that is treated as a population, via re-sampling it, like as in Efron \cite{Efron}, for example, or by re-weighing methods in general.
\par
In Section \ref{sec4}, we demonstrate the validity of the bootstrapped $t$-intervals for both approaches to conditioning. In particular, when conditioning on the sample, our results in this regard are new in that they are shown to hold true when $X\in DAN$, possibly with infinite variance, while the ones with $E_{X} X^2 <\infty$ when conditioning on the weights are first time results \emph{per se}.
 \par
All the proofs are given in Section \ref{sect5}.
\\
\textbf{Notations}. Conditioning on the bootstrap weights $v_i^{(n)}$ and conditioning on the data $X_i$, call for  proper notations that distinguish the two approaches.  Hence, the notation $(\Omega_X,\mathcal{F}_X, P_X)$ will stand for the probability space on which  $X,X_1,X_2,\ldots$ are defined, while $(\Omega_v, \mathcal{F}_v, P_v)$ will stand for the probability space on which the triangular arrays of the  bootstrap weights $v_1^{(1)},(v_1^{(2)}, v_2^{(2)})$, \ldots,$(v_1^{(n)},\ldots,v_n^{(n)}),\ldots$ are defined. In view of the independence of these sets of random variables, jointly they live on the direct product probability space $(\Omega_{X}\times \Omega_{v},\mathcal{F}_{X}\otimes\mathcal{F}_{v},P_{X,v}=P_X\times P_{v})$.    Moreover, for  use throughout, for each $n\geq 1$,  we let $P_{.|v}(.)$ be a short hand notation for the conditional probability $P(.|\mathcal{F}^{(n)}_{v})$ and, similarly, $P_{.|X}(.)$ will stand for the conditional probability $P(.|\mathcal{F}^{(n)}_{X})$, where  $\mathcal{F}^{(n)}_{v}:= \sigma(v^{(n)}_1,\ldots,v^{(n)}_{n})$ and $\mathcal{F}^{(n)}_{X}:=\sigma(X_1,\ldots,X_n)$, respectively, with corresponding conditional expected values $E_{.|v}$ and $E_{.|X}$. In case of Efron's scheme of bootstrap, we will use $w$ instead of $v$ in all these notations whenever convenient.

\section{CLT via conditioning on the bootstrap weights}\label{sect2}

In this section we explore  the asymptotic behavior of the   weighted bootstrap     via conditioning on the bootstrap weights. The major motivation for  conditioning on the weights is that, when bootstrapping the i.i.d.   observables   $X, X_{1},X_{2},\ldots$, these random variables    should continue to be  the prime source of stochastic  variation and, hence,  the random samples should be the main contributors to establishing conditional CLT's for the bootstrapped $t$-statistics as defined in  (\ref{weightedboot1}) and (\ref{weightedboot2}).
The following Theorem \ref{thm1} formulates   the main approach  of this paper to the area of weighted bootstrap for the Student $t$-statistic. Based on a direct consequence of the Lindeberg-Feller CLT (cf. Lemma \ref{new}), it amounts to concluding  appropriate equivalent Lindeberg- Feller type CLT's respectively, corresponding  to both versions of the following statement: as $n,m_n \to \infty$,

$$ M_n:=\frac{ \max_{1\leq i\leq n} \big( \frac{v_i^{(n)}}{m_{n}} - \frac{1}{n}\big)^2 }{\sumn \big( \frac{v_i^{(n)}}{m_{n}} - \frac{1}{n} \big)^2}= \left\{      \begin{array}{ll}
                o(1) \ a.s.-P_v  \\
                o_{P_v}(1).
               \end{array}
\right.
$$

\begin{theorem}\label{thm1}
Let $X, X_{1},X_{2},\ldots$ be real valued i.i.d. random variables with mean $0$ and variance   $\sigma^2,$  and assume that  $0<\sigma^2 < \infty$. Put $V_{i,n}:=  \big| \big( \frac{v^{(n)}_{i}}{m_n}-\frac{1}{n} \big) X_i \big|, \ 1\leq i \leq n$,  $V^{2}_{n}:= \sum_{i=1}^{n}\big( \frac{v^{(n)}_{i}}{m_n}-\frac{1}{n} \big)^{2}$, $M_n:=\frac{ \max_{1\leq i\leq n} \big( \frac{v_i^{(n)}}{m_{n}} - \frac{1}{n}\big)^2 }{\sumn \big( \frac{v_i^{(n)}}{m_{n}} - \frac{1}{n} \big)^2}$,  and let $Z$ be a standard normal random variable throughout. Then, as $n,m_n \to \infty$, having

\begin{eqnarray}
&& M_n=o(1)\ a.s.-P_v   \label{A1}   \\
&&is\ equivalent\ to\ concluding\ the\ respective\ statements\ of\ (\ref{A2})\ and \ (\ref{A3})\nonumber \qquad \qquad \qquad \\ &&simultaneously\ as \ follows \nonumber\\
&& P_{X|v} \left(T^*_{m_n} \leq t\right) \longrightarrow
P(Z\leq t)\ a.s.- P_{v}\ for \ all\  t \in \mathds{R} \label{A2}\\
&&and \nonumber \\
&&\max_{1\leq i \leq n} P_{X|v}(V_{i,n}\big/ (S_{n} V_{n}) >\ep)=o(1)\ a.s.-P_{v}, for \ all\  \ep>0, \label{A3}
\end{eqnarray}
and, in a similar vein, having
\begin{eqnarray}
&& M_n=o_{P_v}(1)  \label{B1}  \\
&& is\ equivalent\ to\ concluding\ the\ respective\ statements\ of\ (\ref{B2})\ and \ (\ref{B3})\ as\ below \qquad     \nonumber\\
&& simultaneously \nonumber\\
&& P_{X|v} \left(T^*_{m_n} \leq t\right) \longrightarrow
P(Z\leq t)\ in\ probability - P_{v}\ for \ all\  t \in \mathds{R} \label{B2}\\
&&and \nonumber\\
&&\max_{1\leq i \leq n} P_{X|v}(V_{i,n}\big/ S_{n} V_{n}>\ep)=o_{P_v}(1),for \ all\  \ep>0.  \label{B3}
\end{eqnarray}
Moreover,    assume that, as $n, m_n \to \infty$, we have for\ any\ $\ep>0$,
\begin{numcases}{
\ P_{X|v}\Big(\big|\frac{S_{m_n}^{*2}\Big/ m_n}{ \sigma^2\ \sumn\big( \frac{v_i^{(n)}}{m_n}-\frac{1}{n} \big)^2} -1 \big|>\varepsilon \Big)=}
o(1)\ a.s.-P_{v} \label{E}& \\
o_{P_{v}}(1). \label{F}&
\end{numcases}
Then, as  $n, m_n \to \infty$,  via (\ref{E}),   the  statement of (\ref{A1}) is also equivalent to having (\ref{c2}) and (\ref{c3}) simultaneously as below

\begin{eqnarray}
&& P_{X|v} \left(T^{**}_{m_n} \leq t\right) \longrightarrow
P(Z\leq t)\ a.s.- P_{v}\ for \ all\  t \in \mathds{R} \label{c2} \qquad \qquad  \qquad \qquad \qquad \qquad \qquad   \\
&& and \nonumber\\
&&\max_{1\leq i \leq n} P_{X|v}(V_{i,n}\big/ (S^{*}_{m_n}/\sqrt{m_n} )>\ep)=o(1)\ a.s.-P_{v}, \ for\ all\ \ep>0, \label{c3}
\end{eqnarray}
and, in a similar  vein, via (\ref{F}), the statement (\ref{B1}) is also equivalent to having  (\ref{d2}) and (\ref{d3}) simultaneously as below

\begin{eqnarray}
&& P_{X|v} \left(T^{**}_{m_n} \leq t\right) \longrightarrow
P(Z\leq t)\ in\ probability - P_{v}\ for \ all\  t \in \mathds{R} \label{d2} \qquad \qquad \qquad \qquad \qquad\\
&&and   \nonumber \\
&&\max_{1\leq i \leq n} P_{X|v}(V_{i,n}\big/ (S^{*}_{m_n}/\sqrt{m_n} )>\ep)=o_{P_v}(1),\ for \ all \ \ep>0.  \label{d3}
\end{eqnarray}
\end{theorem}
For verifying the technical conditions (\ref{E}) and (\ref{F})  as above,     one does not need to know the actual finite   value  of $\sigma^2$.
\par
The essence of Theorem \ref{thm1} is that for i.i.d. data with a finite second moment, a
  scheme of bootstrap for the Student $t$-statistic   is valid \emph{if and only if}  the  random weights in hand  satisfy either one of the maximal negligibility  conditions as in  (\ref{A1}) or (\ref{B1}) for $M_n$. Thus, when conditioning on the weights, Theorem \ref{thm1} provides  an overall approach for obtaining CLT's for bootstrap means in this context, a role that is similar to  that of Theorem 2.1  of Mason and Newton \cite{Mason and Newton} that provides CLT's for generalized bootstrap means of exchangeable weights when conditioning on the sample. Incidentally, conclusion (\ref{B2})  of our Theorem \ref{thm1} under the maximal negligibility conclusion  (\ref{B1}) is a non-parametric version of the scaler scaled (not self-normalized) Theorem 3.1 of Arena-Guti\'{e}rrez  and Mart\'{a}n \cite{Arenal}  under their  more restrictive conditions E1-E5 for exchangeable weights, where condition E4 and E5 combined yield our condition (\ref{B1}) in terms of exchangeable weights. In this regard we also note in passing that, at the end of  Section 1.2 of his lectures on some aspects of the bootstrap \cite{Gine}, Gin\'{e}  notes that checking conditions E4-E5 of \cite{Arenal} sometimes require ingenuity.

\par
When the scheme of bootstrap is specified to be Efron's,  then  Corollary \ref{cor1} hereupon to Theorem \ref{thm1} implies the validity of this scheme for both $T^{*}_{m_n}$ and  $T^{**}_{m_n}$as follows.

\begin{corollary}\label{cor1}
Consider $v^{(n)}_i=w^{(n)}_i$, $1\leq i \leq n$,  $n\geq 1$,  and $M_n$ of Theorem \ref{thm1} in terms of these re-sampling weights as in Remark \ref{remark1}, i.e., Efron's scheme of bootstrap. Assume that $0< \sigma^{2}= var(X) <\infty$.
\\
(a) If  $m_{n}, n \rightarrow \infty$, in such a way that $m_n=o(n^2)$, then, mutatis mutandis, (\ref{B1}) is equivalent to having  (\ref{B2})  and (\ref{B3})  simultaneously,  and spelling out only  (\ref{B2}), in this context  it reads
\begin{equation}
P_{X|w} (T^{*}_{m_n} \leq t) \longrightarrow
  P(Z\leq t)\ in \ probability-P_w \ for\ all\ t \in \mathds{R},
\end{equation}
(b) If  $m_{n}, n \rightarrow \infty$ in such a way that $m_n=o(n^2)$ and $ n=o(m_n)$, then, mutatis mutandis again,  (\ref{B1}) is also equivalent to having
(\ref{d2}) and (\ref{d3}) simultaneously, and spelling out only    (\ref{d2}), in this context it reads as follows
\begin{equation}
 P_{X|w} (T^{**}_{m_n} \leq t) \longrightarrow
  P(Z\leq t)\ in \ probability-P_w, for\ all\ t \in  \mathds{R}.
 \end{equation}
\end{corollary}

\begin{remark}\label{S_n fixed}
It is noteworthy to note that, along the lines of the proof of the preceding corollary (cf. the second part of the proof of our Lemma \ref{lemma3}), it will be seen that for a finite number of observations $X_1,\ldots,X_n$ in hand, $S^{*^2}_{m_n}$, i.e., the bootstrap version of the sample variance $S_{n}^{2}$,  is an in probability-$P_{X,w}$ consistent   estimator of $S^{2}_{n}$, as only $m_n=: m \to \infty$.  In other words, when $E X^{2}_{1}<\infty$,  on taking $n$ to be fixed as $m_n=m \to \infty$, we have that
\begin{equation}\label{S_n fixed consisitency}
S^{*^2}_{m_n} \longrightarrow S_{n}^{2}\ in \ probability - P_{X,w}.
\end{equation}
Consequently, the bootstrap sample variance of only one large enough bootstrap sub-sample yields a consistent estimator for  the sample variance $S_{n}^{2}$ of the original sample. Moreover, a similar result can be shown to also  hold true for estimating the mean of the original sample  $\bar{X}_{n}$ via taking only one large enough bootstrap sub-sample and computing its mean $\bar{X}_{m_n}^{*}$ when $n$ is fixed. In fact, in a more general setup,  the consistency  result (\ref{S_n fixed consisitency})   for characteristics of the original sample which are of the form of $U$-statistics can be found in Cs\"{o}rg\H{o} and Nasari \cite{Csorgo and Nasari} (cf. Part (a) of Theorem 3.2). These results provide an alternative   to   the classical method, as suggested, for example,  by Efron and Tibshirani \cite{EfronTibshirani}, where  the average of the bootstrapped estimators, $\bar{X^*}(b)$ of $B$ bootstrap sub-samples drawn repeatedly and independently   from the original sample, is considered as an estimator for a characteristic  of the sample in hand,  such as $\bar{X}_n$  and $S^{2}_n$, for example. The  validity of the average of these $B$ bootstrap estimators is then investigated as  $B\to \infty$.
\end{remark}

\begin{remark}\label{remark8}
In probability-$P_w$, part (b) of Corollary \ref{cor1}   parallels (1.11) of  Theorem 1.1 of Mason and Shao \cite{MasonShao} in which they conclude  that, when $E X^2 <\infty$, then  for {\em almost all realizations} of the sample (i.e., for almost all samples), the conditional (on the data) distribution of  $T^{**}_{m_n}$   will coincide with the standard normal distribution whenever $\lambda_1\leq m_n/ n \leq \lambda_2$ for all $n$ large enough and some constants $0< \lambda_1 < \lambda_2< \infty$. It would be desirable to have an a.s.-$P_w$ version of our Corollary \ref{cor1}, and to extend the in probability$- P_w$  validity of its present form to having  $X\in DAN$ with $E X^2 =\infty$.
\end{remark}
\par
Now suppose that  $\displaystyle{v^{(n)}_{i}=\zeta_{i} }$, $1\leq i \leq n$, where $\zeta_{i}$ are positive i.i.d. random variables. In this case the bootstrapped $t$-statistic $T^{*}_{m_{n}}$ defined by (\ref{weightedboot1}) is of the form:
\begin{equation}
T^*_{m_{n}} =
\frac{ \displaystyle \sum^n_{i=1} \big( \frac{\zeta_{i}}{m_n} - \frac{1}{n} \big) X_i}
{S_n \sqrt{   \displaystyle\mathop{\sum}_{i=1}^{n} (\frac{\zeta_{i}}{m_n}-\frac{1}{n})^2 } },\label{T*zeta}
\end{equation}
where $\displaystyle{m_n=\sum_{i=1}^{n} \zeta_i}$.
\par
 The following Corollary \ref{cor2} to  Theorem \ref{thm1}  establishes the validity of this scheme of bootstrap for $T^{*}_{m_{n}}$, as defined by (\ref{T*zeta}), via conditioning on the  bootstrap weights of the latter.
\begin{corollary}\label{cor2}
Assume that $0 < \sigma^{2}= var(X) <\infty$, and let $\zeta_1,\zeta_2,\ldots$ be a sequence of positive i.i.d. random variables  which are independent of $X_{1},X_{2}, \dots$ .  Then,    as $ n\to\infty$,
\\
(a) if   $E_{\zeta}(\zeta^{4}_{1})< \infty$,  then, mutatis mutandis, condition (\ref{A1}) is equivalent to having (\ref{A2}) and (\ref{A3}) simultaneously, and spelling  out only  (\ref{A2}),   in this context it reads
\begin{equation}
P_{X|\zeta} (T_{m_{n}}^* \leq t) \longrightarrow P(Z\leq t) \ \hbox {a.s.}-P_\zeta, \ for\ all\ t \in \mathds{R},
\end{equation}
(b) if   $E_{\zeta}(\zeta^{2}_{1})< \infty$,  then, mutatis mutandis, (\ref{B1}) is equivalent (\ref{B2}) and (\ref{B3})  simultaneously, and spelling  out only  (\ref{B2}),   in this context it reads
\begin{equation}
P_{X|\zeta} (T_{m_{n}}^* \leq t) \longrightarrow P(Z\leq t) \ in \ probability-P_\zeta, \ for\ all\ t \in \mathds{R},
\end{equation}
where $Z$ is a standard normal random variable.
\end{corollary}
\section{CLT via conditioning on the sample}\label{sect3}
Efron's bootstrapped partial sums via conditioning on the data have been the subject of intensive study and   many remarkable papers  can be found in the literature in this regard.
\par
Conditioning on the data which are assumed to be in DAN, Hall \cite{Hall} proved that if $m_n, \, n\to\infty$, and  $\lambda_1 \leq m_n\big/n \leq \lambda_2$, where $0< \lambda_1 < \lambda_2 < \infty$, then there exists a sequence of positive numbers $\left\{\gamma_n\right\}^\infty_{n=1}$ such that
\begin{equation}\label{eq4}
\frac{\sqrt{m_n}(\bar{X}^*_{m_n} - \bar X_n)}{\gamma_n}  \dto N(0,1) \hbox{ ~in~ probability}- P_X.
\end{equation}
In the same year S. Cs\"{o}rg\H{o} and Mason \cite{CsorgoMason} showed that under the same conditions as those assumed by Hall, i.e., $X\in DAN$ and $m_n\big/ n \in [\lambda_1,\lambda_2]$ with $0<\lambda_1< \lambda_2<\infty$ as before, the numerical constants $\gamma_n$ in (\ref{eq4})
can be replaced by the sample standard deviation $S_n$, and the conclusion of   (\ref{eq4}) remains true.   Furthermore, Mason and Shao \cite{MasonShao} replaced $S_n$ by the bootstrapped sample standard deviation $S^*_{m_n}$ and, under the  conditions assumed by Hall \cite{Hall} and  S. Cs\"{o}rg\H{o} and Mason \cite{CsorgoMason}, i.e., when $m_n/n \in [\lambda_1,\lambda_2]$, they concluded  that
\begin{equation}\label{eq5}
T^{**}_{m_n} \dto N(0,1) \hbox{ in ~ probability}-P_{X}
\end{equation}
 if and only if $X\in DAN$, possibly  with $E X^2 =\infty$. As mentioned already (cf. Remark \ref{remark8}),  when $m_n/n \in [\lambda_1,\lambda_2]$,  Mason and Shao \cite{MasonShao} also  characterized  the almost sure-$P_X$ validity (asymptotic normality) of $T^{**}_{m_n}$ via conditioning on the data when their variance is positive and finite.
 \par
 Thus, whenever $m_n/n \in [\lambda_1,\lambda_2]$, via conditioning on the data which are in $DAN$,  Mason and Shao \cite{MasonShao} established the validity in probability-$P_X$ of the Efron bootstrapped version of the $t$-statistics as in (\ref{weightedboot2}), as well as its almost sure$-P_X$  validity when $E X^2$ is  positive and finite (cf. (1.10) and (1.11), respectively,  of their Theorem 1.1). Under its condition (\ref{Newly added}) the respective conclusions of our (\ref{eqadded1}) and (\ref{eqadded1prime}) of our forthcoming Theorem \ref{theorem2} parallel those of (1.10) and (1.11) of Theorem 1.1 of Mason and Shao \cite{MasonShao},   who  also noted the desirability of having (\ref{eq5}) holding true when the data are in $DAN$ and $m_n=n$.  Theorem \ref{theorem3} below relates to this question in terms of $(S_{m_n}^{*}\big/ S_{n}) T_{m_n}^{**}$ (cf. (\ref{added(3.6)}) and  Remark \ref{remark9}).
\begin{remark}\label{remark5}
 For a rich source of information  on the topic of bootstrap we refer to  the insightful  survey by S. Cs\"{o}rg\H{o} and Rosalsky \cite{CsorgoRosalsky}, in which various types of limit laws are  studied for bootstrapped sums.
\end{remark}
 Among those who explored   weighted bootstrapped partial sums, we mention S. Cs\"{o}rg\H{o} \cite{Csorgo} and Arenal-Guti\'{e}rrez \emph{et al}. \cite{Arenal}, who studied the  unconditional strong law of large numbers for the bootstrap mean.

 Mason and Newton \cite{Mason and Newton} introduced the idea of the generalized bootstrap for the sample mean that is to replace the multinomial Efron bootstrap as in our Remark \ref{remark1} by another vector of exchangeable non-negative random variables that are also independent of the $X_i$. Their basic  tool for establishing the almost sure-$P_X$ CLT consistency  of their generalized bootstrap mean, as in their Theorem 2.1, is Theorem 4.1 of H\'{a}jek   \cite{Hajeck} concerning the asymptotic normality of linear rank statistics. Accordingly,  their Theorem 2.1 deals with the a.s.-$P_X$ asymptotic normality of  exchangeable arrays of self-normalized partial sums when conditioning on the sample.

 %a new method to study convergence in distribution  of bootstrapped partial  sums by first  deriving a CLT for them in terms of the joint distribution $P_{X,w}$. Their approach relies on a result due to Hajeck \cite{Hajeck} concerning exchangeable random variables and it is a parallel result to Theorem \ref{thm1} in the present exposition.
 \par
 Taking a different approach  form that of Mason and Newton \cite{Mason and Newton}, Arenal-Guti\'{e}rrez and Matr\'{a}n  \cite{Matran} developed a technique by which they derived  a scaler scaled almost sure-$P_X$, conditional on the sample, CLT for $(\bar{X}_{m_{n}}^{*}-\bar{X}_n)$ with the parametric scaler $\sqrt{var(X)}/\sqrt{m_n}$ (cf. their Theorem 3.2)

 \par
  Conditioning on the sample, in this section we study the validity of Efron's scheme of bootstrap when applied to sums of i.i.d. random variables. As will be seen, in establishing a conditional CLT, given the data,   the weights, $w_{i}^{(n)}$, as random variables, weighted by conditioning on the data,  will play the dominant role. This is in contrast to the previous section,  in which  a weighted i.i.d. version of the Lindeberg-Feller CLT for the data, $X$, played the dominant role in deducing our Theorem \ref{thm1}.
  \par
Clearly (cf., e.g., Lemma 1.2 in S. Cs\"{o}rg\H{o}  and Rosalsky \cite{CsorgoRosalsky}),  \emph{unconditional central limit theorems} result from the conditional ones in $P_v$ or $P_X$ under their respective conditions,  and, in turn,  this is the way bootstrap works when taking repeated bootstrap samples (cf. our Section \ref{sec4}). S. Cs\"{o}rg\H{o} and Rosalsky \cite{CsorgoRosalsky} indicate that the laws of unconditional bootstrap are ``less frequently spelled out in the literature".  Hall \cite{Hallunconditional}, however, addresses both conditional and unconditional laws for bootstrap.   S. Cs\"{o}rg\H{o} and Rosalsky \cite{CsorgoRosalsky} also note that, according to Hall,  conditional laws are of interest to statisticians who are interested in the probabilistic aspects of the sample in hand, while the unconditional laws of bootstrap have the ``classical frequency interpretation". Accordingly, and as noted already, our approach in Section \ref{sect2} is that of a statistician interested in studying the probabilistic aspects  of a sample that is treated as a population, by means of conditioning on re-sampling, and/or, re-weighing the data in hand.
\par
We wish  to emphasize  that in this section \emph{only} Efron's scheme of bootstrap will be considered.  This is so, since the validity and establishment of the results here, to a large extent,  rely on the multinomial structure of the random weights, $w^{(n)}_{i}$, in this scheme.  On the other hand, the data are  assumed to be in $DAN$, possibly with infinite variance,  and studied   under conditions on $n, m_n$, as $n\rightarrow \infty$, that differ from requiring $m_n/n$ to be  in the interval $[\lambda_1,\lambda_2]$ with $0< \lambda_1<\lambda_2<\infty$ as in Mason and Shao \cite{MasonShao}.
\par
It is well-known that the $t$-statistic converges in distribution  to a standard normal random variable   if and only if the data are in $DAN$ (cf. Gin\'{e} \emph{et al.}
\cite{GineGotzeMason}).
The following Theorem \ref{theorem2} establishes the validity (asymptotic normality) of the Efron bootstrapped   version of the $t$-statistics as in (\ref{weightedboot2}), based on random samples on $X\in DAN$ via conditioning on the data. It is to be compared to the similarly conditioned  Theorem 1.1 of Mason and Shao \cite{MasonShao}.

\begin{theorem}\label{theorem2}
Let $X, X_1,\ldots$ be i.i.d.  random variables with $X\in  DAN$.
Consider $T^{**}_{m_n}$ as in (\ref{weightedboot2}) with $X \in DAN$ and   Efron's bootstrap $\{ w^{(n)}_i, 1\leq i \leq n\}$, ${n\geq 1}$, scheme of re-sampling from random samples $\{ X_i, 1\leq i \leq n \}_{n\geq 1}$ as in (\ref{correction1}) and Remark \ref{remark1}.  If, as $n, m_n\to\infty$ so that
\begin{equation} \label{Newly added}
 \frac{m_n}{ 2 n \log n }  \to  \infty,
\end{equation}
then, $\ for\ all \ t\in \mathds{R}$,
\begin{equation}\label{eqadded1}
P_{w|X} \big(T^{**}_{m_n} \leq t\big) \longrightarrow P(Z\leq t) \hbox{ in \ probability}-P_X,
\end{equation}
and, when $E_{X} X^2< \infty$, then
\begin{equation}\label{eqadded1prime}
P_{w|X} \big(T^{**}_{n} \leq t\big) \longrightarrow P(Z\leq t) \ a.s.-P_{X},
\end{equation}
where, $Z$ is a standard normal random variable. Further to (\ref{eqadded1prime}), if $n,m_n \rightarrow \infty$ so that, instead of (\ref{Newly added}), we have
\begin{equation}\label{Newly added2}
m_n/n \rightarrow \infty
\end{equation}
then, when $E_{X} X^2<\infty$, (\ref{eqadded1prime}) continues to hold true in probability-$P_X$.
\end{theorem}
\par
\par
The next result relates to a question raised by Mason and Shao \cite{MasonShao} asking if  the conditional CLT in (\ref{eqadded1}) held
 true when $m_n=n$.  According to the following Theorem \ref{theorem3}, the answer is positive if one replaces  $T^{**}_{m_n}$ by
 \begin{equation}
 \label{eqTSn}
 T^{**}_{m_n,S_{n}}:= \frac{\sum_{i=1}^{n} (\frac{w_{i}^{(n)}}{m_n}-\frac{1}{n} ) X_{i} }{S_{n}\big/ \sqrt{m_n}}= \frac{S^{*}_{m_n}}{S_n} T^{**}_{m_n}.
 \end{equation}
\begin{theorem}\label{theorem3}
Let $X, X_1,\ldots$ be  i.i.d. random variables with $X\in DAN$.
Consider Efron's bootstrap scheme as in Theorem \ref{theorem2}.   If, as $n, \ m_n\to\infty$, so that  for an  arbitrary $\ep>0$  we have $\frac{m_n}{n}\geq \ep >0 $,
then, for all $t\in \mathds{R}$,
\begin{equation}\label{added(3.6)}
P_{w|X} \big(T^{**}_{m_n,S_{n}} \leq t\big) \longrightarrow P(Z\leq t) \hbox{ in \ probability }-P_X,
\end{equation}
where  $Z$ is a standard normal random variable.
\end{theorem}
\begin{remark}\label{remark9}
On  taking $m_n=n$, Theorem \ref{theorem3} continues to hold true as before, but now  in terms of
\begin{equation}
T^{**}_{n,S_{n}} = \frac{S^*_{n}}{S_n} \ T^{**}_{n}.\nonumber
\end{equation}
\end{remark}

\par
The conclusion of (\ref{added(3.6)}) coincides  with that of (5.2) of S. Cs\"{o}rg\H{o} and Mason \cite{CsorgoMason}, who,  as mentioned right after (\ref{eq4}) above, concluded it for $X \in DAN$ whenever, $m_n/n \in [\lambda_1,\lambda_2]$ with $0<\lambda_1<\lambda_2<\infty$. Thus, for our conclusion in (\ref{added(3.6)}), we may take $m_n/n \in [\lambda_1,\infty)$, and conclude also Remark \ref{remark9} with $m_n=n$ that was first established by Athreya \cite{Athreya}. For further comments  along these lines we refer to Section 5 of S. Cs\"{o}rg\H{o} and Mason \cite{CsorgoMason}.

\section{Validity  of Bootstrapped $t$-intervals}\label{sec4}
In order  to establish an asymptotic confidence bound for $\mu=E(X)$ with an asymptotic probability coverage  of size $\alpha$, $0< \alpha\leq 1$, using the classical CLT, one can use the classical Student pivot $T_{n}$ via setting
 $T_{n}\leq z_{\alpha}$, where $P(Z\leq z_{\alpha} )=\alpha$. One can also establish an asymptotic  size $\alpha$ bootstrap confidence bound for $\mu$  by taking $B\geq 1$ bootstrap sub-samples of size $m_{n}$ via re-sampling,    or by generating $B$ sets of stochastically reweighed bootstrap sub-samples of $\{X_i, 1\leq i \leq n \}$ independently (i.e., each set of the $B$ bootstrap weights are independent). The latter can be done by simulating  $B$ sets of independent i.i.d. weights $(\zeta_{1}^{(b)},\dots,\zeta_{n}^{(b)})$, $1\leq b \leq B$.    Obviously, the independence of the bootstrap weights with respect to the probability $P_{v}$ does not imply the independence of  the thus generated  sub-samples with respect to the joint distribution of the data and the bootstrap weights. One will have $B$ values of $T_{m_{n}}^*(b)$  and/or $T_{m_{n}}^{**}(b)$ or $T_{n}^{*}(b)$, $1\leq b \leq B$,  and respective       asymptotic $100.\alpha\% $ bootstrap confidence bounds  will result, as in the upcoming Theorems \ref{Efron confidence}  and \ref{zeta confidence}, from the inequalities
\begin{equation}\label{bootstrapconfidencebound}
T_{n}\leq C^{(B)}_{s,\alpha},\ s=1,2,3,4
\end{equation}
$ \vspace{-.5 cm}$
 where
\begin{eqnarray*}
C^{(B)}_{1,\alpha}&:=&\inf\{t: \ \frac{1}{B}\sum_{ b=1}^{B} I(T_{m_{n}}^{*}(b)\leq t)\geq \alpha \}, \\
C^{(B)}_{2,\alpha}&:=&\inf\{t: \ \frac{1}{B}\sum_{ b=1}^{B} I(T_{m_{n}}^{**}(b)\leq t) \geq \alpha \},\\
C^{(B)}_{3,\alpha}&:=&\inf\{t: \ \frac{1}{B}\sum_{ b=1}^{B} I(T_{m_{n},S_{n}}^{**}(b)\leq t) \geq \alpha \},  \\
C^{(B)}_{4,\alpha}&:=&\inf\{t: \ \frac{1}{B}\sum_{ b=1}^{B} I(T_{n}^{*}(b)\leq t)\geq \alpha \},
\end{eqnarray*}
and $T_n$ is the Student $t$-statistic as in (\ref{t-statistic}).
\par
Observe that $C^{(B)}_{s,\alpha},\ s=1,2,3,4$,  are  bootstrap estimations of the respective  $100.\alpha$ percentile of the  distributions $P_{X,v} \big(T^{*}_{m_n} \leq t\big) $, $P_{X,v} \big(T^{**}_{m_n} \leq t \big)$, $P_{X,v} \big(T^{**}_{m_n,S_n} \leq t \big)$ and $P_{X,v} \big(T^{*}_{n} \leq t \big)$. Moreover, since $\displaystyle{C^{(B)}_{s,\alpha}}$ are the $100.\alpha$ percentiles  of their respective  empirical distributions, therefore they  coincide  with  their respective  order statistics $T^{*^{(l)}}_{m_{n}}$,  $T^{**^{(l)}}_{m_{n}}$, $T^{**^{(l)}}_{m_{n},S_n}$ and   $T^{*^{(l)}}_{n}$,  where $l=[\alpha(B+1)]$.  \par
We note that  $\displaystyle{C^{(B)}_{s,\alpha},\ s=1,2,3,4}$,  are  natural extensions of S. Cs\"{o}rg\H{o} and Mason's \cite{CsorgoMason} approach to  establishing the validity of bootstrapped empirical processes. Some ideas  that are used in the proofs of the results in this section were borrowed from  \cite{CsorgoMason} and adapted accordingly.
 \par
 The objective of this section is to show that in the light of Theorems \ref{thm1}, \ref{theorem2} and \ref{theorem3},    the confidence bounds  obtained from  (\ref{bootstrapconfidencebound}) will achieve  the nominal coverage probability  $\alpha$ as $n, \ m_{n}$  and $B\rightarrow \infty$. More precisely, in  Theorem  \ref{Efron confidence} below we consider the confidence bound as in  (\ref{bootstrapconfidencebound}) and   Efron's scheme of bootstrap,  and   show that  the asymptotic nominal coverage probability  $\alpha$ will be achieved. Moreover, the latter will be shown to be true via  conditioning on the bootstrap weights and also via conditioning on the data. In Theorem \ref{zeta confidence} we consider the confidence bound in (\ref{bootstrapconfidencebound}) with $C^{(B)}_{4,\alpha}$  when the scheme of bootstrap is stochastically re-weighing and via conditioning on  the bootstrap weights, we show that the  asymptotic nominal coverage probability $\alpha$ will again  be achieved.

\par
Thus, both approaches to the bootstrap will be shown  to work, namely, as in (a) of Theorem \ref{Efron confidence} and as in Theorem \ref{zeta confidence} when conditioning on the weights, and as in (b) and (c) of Theorem \ref{Efron confidence} when conditioning on the data.

\par
In order to state the just mentioned  conclusions, one needs to define an appropriate  probability space for accommodating the  presence of $B$ bootstrap sub-samples, as $B\to \infty$. This means that one has to incorporate  $B$ i.i.d. sets of weights $$\Big(v_{1}^{(1)}(b),\textbf{(}v_{1}^{(2)}(b),v_{2}^{(2)}(b)\textbf{)},\ldots,\textbf{(}v_{1}^{(n)}(b),\ldots,v_{n}^{(n)}(b)\textbf{)}
,\ldots  \Big),$$
which live on their respective   probability spaces $(\Omega_{v(b)},\mathfrak{F}_{v(b)}, P_{v(b)})$, $b\geq 1$.
In view of  this, and due to the fact that $n,\ m_{n}$ and  $B$ will approach $\infty$, we let $(\bigotimes_{b=1}^{\infty}\Omega_{v(b)}, \bigotimes_{b=1}^{\infty}\mathfrak{F}_{v(b)}$ $,\bigotimes_{b=1}^{\infty}P_{v(b)})$   be the probability space on which the following  row-wise i.i.d. array of bootstrap weights are  defined:
$$
\begin{array}{cccc}
v_{1}^{(1)}(1),&\textbf{(} v_{1}^{(2)}(1),v_{2}^{(2)}(1)\textbf{)},& \textbf{(} v_{1}^{(3)}(1),v_{2}^{(3)}(1), v_{3}^{(3)}(1)\textbf{)}  ,&\ldots \\
v_{1}^{(1)}(2),&\textbf{(} v_{1}^{(2)}(2),v_{2}^{(2)}(2)\textbf{)},& \textbf{(} v_{1}^{(3)}(2),v_{2}^{(3)}(2), v_{3}^{(3)}(2)\textbf{)}  ,&\ldots \\
 \vdots & \vdots & \vdots & \vdots
\end{array}
$$
In what  follows,  we let $ (\bigotimes_{b=1}^{\infty}\Omega_{X,v(b)}, \bigotimes_{b=1}^{\infty}\mathfrak{F}_{X,v(b)},\bigotimes_{b=1}^{\infty}P_{X,v(b)})  $ be the joint probability space of the $X$'s and the preceding array of the  weights $v(b)$, $b\geq 1$.
\begin{thm}\label{Efron confidence}
 Consider Efron's scheme of bootstrap, i.e., $v_{i}^{(n)}=w_{i}^{(n)}$, $1\leq i\leq n$, $n\geq 1$.
 \\

(a) Assume the  conditions of Corollary \ref{cor1}. Then, as $n,m_n,B\rightarrow \infty$,
$$C^{(B)}_{1,\alpha},C^{(B)}_{2,\alpha}\longrightarrow z_{\alpha}\ in\ probability  - \bigotimes_{b=1}^{\infty}P_{X,w(b)} .$$

(b) Assume the  conditions of Theorems \ref{theorem2}. Then,  as $n,m_n,B\rightarrow \infty$,
$$ C^{(B)}_{2,\alpha}\longrightarrow z_{\alpha}\ in\ probability -  \bigotimes_{b=1}^{\infty}P_{X,w(b)}.
$$

(c) Assume the  conditions of Theorem \ref{theorem3}. Then,  as $n,m_n,B\rightarrow \infty$,
$$ C^{(B)}_{3,\alpha}\longrightarrow z_{\alpha}\ in\ probability -  \bigotimes_{b=1}^{\infty}P_{X,w(b)}.
$$
\end{thm}

\begin{thm}\label{zeta confidence}
Suppose that $v_{i}^{(n)}=\zeta_i$, $1\leq i \leq n$, and put $m_n=\sumn \zeta_i$. Assume the  conditions  of Corollary \ref{cor2}. Then, as $n,B\rightarrow \infty$,
$$
 C^{(B)}_{4,\alpha}\longrightarrow z_{\alpha}\ in\ probability - \bigotimes_{b=1}^{\infty}P_{X,\zeta(b)}.
$$
\end{thm}

\par
When conditioning on the sample, the validity of the bootstrap confidence intervals  was also studied by Hall \cite{Hallunconditional} when, with some $\delta>0$, $E_{X} X^{4+\delta}<\infty$ and $m_n=n$. Our conclusions in  (b)  and (c) hold true when $X \in DAN$, possibly with infinite variance. Conclusion (a) of Theorem \ref{Efron confidence}  and that of Theorem \ref{zeta confidence} are first time results for establishing the validity of bootstrap confidence intervals via conditioning on the weights when $E_{X} X^2 < \infty$.

\section{Proofs}\label{sect5}

The proof of Theorem \ref{thm1} is based on the following Lemma \ref{new} that amounts to a realization of the Lindeberg-Feller CLT.

\begin{lemma}\label{new}

Let $X,X_1,\ldots$ be real valued i.i.d. random variables with mean $0$ and variance $0<\sigma^2 <\infty$ on $(\Omega_X , \mathfrak{F}_{X}, P_{X})$, as before, and let $\left\{a_{i,n}\right\}^n_{i=1}$, $n\geq 1$ be a  triangular array of real valued constants. Then, as $n\to \infty$, %the following statements are equivalent:
$\vspace{-.4 cm}$
\begin{equation}\label{new0}
M_n= \frac{\max_{1\leq i \leq n} a^{2}_{i,n}  }{\sumn a^{2}_{i,n} }\longrightarrow 0,
\end{equation}
if and only if
\begin{equation}\label{new1}
\frac{\sumn a_{i,n} X_{i}}{\sigma \sqrt{\sumn a^{2}_{i,n} }} \to_{d} N(0,1), and, \ for\ all\ \ep>0,
\ \max_{1\leq i \leq n}P_{X}\big( \frac{|a_{i,n} X_{i}  |}{\sigma \sqrt{\sumn a^{2}_{i,n} }}>\ep \big) \to 0
\end{equation}
or, equivalently, if and only if
\begin{equation}\label{new2}
\frac{\sumn a_{i,n} X_{i}}{S_n \sqrt{\sumn a^{2}_{i,n} }} \to_{d} N(0,1),\ and,\ for\ all \ \ep>0,\
\max_{1\leq i \leq n}P_{X}\big( \frac{|a_{i,n} X_{i}  |}{S_n \sqrt{\sumn a^{2}_{i,n} }}>\ep \big) \to 0,
\end{equation}
where $N(0,1)$ stands for  a standard normal random variable, and $S_n$ is the sample variance of the first $n\geq 1$ of the mean $0$ and variance $\sigma^2$ i.i.d. sequence  $X,X_1,X_2,\ldots$ of random variables.
\end{lemma}
\subsection*{Proof of Lemma \ref{new}}
The equivalence  of the respective two statements of  (\ref{new1}) and (\ref{new2}) is an immediate consequence of Slutsky's theorem via  having $S_{n}^{2} \to \sigma^2$ in probability as $n\to \infty$. Hence, it suffices to establish the equivalence of the statement (\ref{new0}) to  the two simultaneous statements of (\ref{new1}).
\par
First assume that we have (\ref{new0}) and show that it implies Lindeberg's conditions that in our context reads as follows: with $F(x)=P_X (X\leq x)$,
\begin{equation}\label{new3}
L_{n}(\ep):= \frac{1}{\sigma^2 \sumn a^{2}_{i,n} }  \sumn a^{2}_{i,n} \int_{(|a_{i,n} x|>\ep \sigma \sqrt{\sumn a^{2}_{i,n}})} x^2 dF(x) \to 0
\end{equation}
 for each $\ep>0$, as $n\to \infty$. Now observe that $L_{n}(\ep)$ can be bounded above by
 \begin{equation}\label{new4}
\frac{1}{\sigma^2} \int_{(|x|>\ep \sigma \sqrt{\frac{\sumn a^{2}_{i,n}}{\max_{1\leq i \leq n} a_{i,n}^{2}}  }) } x^{2} dF(x) \to 0, \ as \ n\to \infty,
\end{equation}
on assuming (\ref{new0}) and $E X^2 = \int x^2 dF(x)$, i.e., (\ref{new0}) implies (\ref{new3}). The latter, in turn, implies the Lindeberg CLT statement of (\ref{new1}). Moreover, by Chebeshev's inequality, via (\ref{new0}) we conclude also the second, the so-called uniform asymptotic uniform negligibility condition  statement of (\ref{new1}). Thus, we now have that (\ref{new0}) implies (\ref{new1}).
\par
Conversely, on assuming now (\ref{new1}), its Lindeberg-Feller type simultaneous conclusions imply the Lindeberg condition of (\ref{new3}), as per the Lindeberg-Feller CLT, and (\ref{new3}) yields (\ref{new0}). $\square$

\subsection*{Proof of Theorem \ref{thm1}}
In view of Lemma \ref{new}, the a.s.-$P_{v}$ equivalence  of (\ref{A1}) to (\ref{A2})-(\ref{A3}) and, via (\ref{E}), that  of  (\ref{A1}) to  (\ref{c2})-(\ref{c3}) hold true along a set $N\in \mathfrak{F}_{v}$ with $P_{v}(N)=1$.
\par
As for the in probability-$P_{v}$ equivalence   of (\ref{B1}) to (\ref{B2})-(\ref{B3}) and, via (\ref{F}), also to  (\ref{d2})-(\ref{d3}), they hold true via the  characterization of convergence in probability in terms of a.s. convergence of subsequences. Accordingly,   for each subsequence $\left\{ n_k  \right\}_{k}$ of $n$, $n\geq 1$,  there exists  a further subsequence $\left\{ n_{k_{\ell}}\right\}_{\ell}$ along which, as $\ell \to \infty$, by virtue of Lemma \ref{new}, the latter two in probability-$P_{v}$ equivalencies reduce to appropriate a.s.-$P_v$ equivalences. This also  completes the proof of Theorem \ref{thm1}. $\square$

\subsection*{Proof of Corollary \ref{cor1}}
Here the bootstrap weights $v_i^{(n)}=w_i^{(n)}$, $1\leq i \leq n$, $n\geq 1$,  are as in Remark \ref{remark1}, i.e., for each $n\geq 1$,
\begin{equation*}
\big(w_1^{(n)},\ldots,w_n^{(n)}\big)\ \substack{d\\=} \ \hbox{multinomial}\big(m_n,\frac{1}{n},\ldots,\frac{1}{n} \big),
\end{equation*}
with $m_n= \sumn w_i^{(n)}  $. In view of Theorem \ref{thm1},  part (a) of Corollary \ref{cor1} will follow  from the  following Lemma \ref{lemma2}, and Lemmas \ref{lemma2} and \ref{lemma3} together will conclude  part (b). $\square$

We now state and prove Lemmas \ref{lemma2} and Lemma \ref{lemma3}.
\begin{lemma}\label{lemma2}
Consider Efron's scheme of bootstrap and assume that $\sigma^2=var(X)<\infty$. If $m_n,n\rightarrow \infty$ in such a way that $m_{n}=o(n^2 )$, then,
\begin{equation*}
\displaystyle{M_n=\frac{\max_{1\leq i\leq n} \big( \frac{w_i^{(n)}}{m_n} - \frac{1}{n} \big)^2}
       {\sum_{i=1}^{n} \big( \frac{w_i^{(n)}}{m_n} - \frac{1}{n} \big)^2} }
\longrightarrow 0 ~\hbox{ in\ probability}-P_w.
\end{equation*}

\end{lemma}

\begin{lemma}\label{lemma3}
Consider Efron's scheme of bootstrap and assume that $0< \sigma^{2}= var(X) <\infty$.  As $m_n,\, n\to\infty$ in such a way that $m_n=o(n^2)$ and $n=o(m_n)$, then,
$$
 P_{X|w}\Big(\big|\frac{S_{m_n}^{*2}\Big/ m_n}{ \sigma^2\ \sumn\big( \frac{w_i^{(n)}}{m_n}-\frac{1}{n} \big)^2} -1 \big|>\varepsilon \Big)\to 0 ~\hbox{ in\ probability}-P_w.
$$
\end{lemma}
\subsection*{Proof of Lemma \ref{lemma2}}
In order to prove this lemma,  for $\varepsilon, \varepsilon^{\prime} > 0$, we write:
\begin{eqnarray*}
&&P_{w}\big(    \frac{\max_{1\leq i\leq n} \big( \frac{w_i^{(n)}}{m_n} - \frac{1}{n} \big)^2}
       {\sum_{i=1}^{n} \big( \frac{w_i^{(n)}}{m_n} - \frac{1}{n} \big)^2}>\ep  \big)\\
&\leq&    P_{w}\big(    \frac{\max_{1\leq i\leq n} \big( \frac{w_i^{(n)}}{m_n} - \frac{1}{n} \big)^2}
       {\sum_{i=1}^{n} \big( \frac{w_i^{(n)}}{m_n} - \frac{1}{n} \big)^2}>\ep, \big| \frac{m_n}{(1-\frac{1}{n})} \sum_{i=1}^{n} \big( \frac{w_i^{(n)}}{m_n} - \frac{1}{n} \big)^2 -1 \big|\leq {\ep}^{\prime}  \big)\\
&+& P_{w}\big(    \big| \frac{m_n}{(1-\frac{1}{n})} \sum_{i=1}^{n} \big( \frac{w_i^{(n)}}{m_n} - \frac{1}{n} \big)^2 -1 \big|> {\ep}^{\prime}  \big)\\
&=& P_w \big( \max_{1\leq i\leq n} \big( \frac{w_i^{(n)}}{m_n} - \frac{1}{n} \big)^2>\frac{\ep(1-\ep^{\prime})(1-\frac{1}{n})}{m_{n}}   \big)
\\
&+& P_w \big(    \big|  \sum_{i=1}^{n} \big( \frac{w_i^{(n)}}{m_n} - \frac{1}{n} \big)^2 -\frac{(1-\frac{1}{n})}{m_n}\big|> \frac{{\ep}^{\prime} (1-\frac{1}{n})   }{m_n}    \big)
\\
&=:& L_1(n)+L_2(n).
\end{eqnarray*}
An upper bound for $L_1(n)$ is:
\begin{eqnarray*}
L_1(n)&\leq&  n   P_w  \big( \big| \frac{w_1^{(n)}}{m_n} - \frac{1}{n}  \big|>\sqrt{\frac{\ep(1-\ep^{\prime})(1-\frac{1}{n})}{m_{n}} }   \big)
\\
&\leq& n \exp\{- \sqrt{m_n}\ . \ \frac{\ep(1-\ep^{\prime})  (1-\frac{1}{n})  }{2\big( \frac{\sqrt{m_n}}{n}+ \sqrt{\ep(1-\ep^{\prime})  (1-\frac{1}{n})  } \big)}  \}.
\end{eqnarray*}
The preceding relation, which  is due to Bernstien's inequality, is a general term of a finite series when $m_n=O(n^2)$.
\par
As for $L_2(n)$, we first note that for each $i$,  $1\leq i \leq n$,
$$\displaystyle{ E_w (\frac{w_i^{(n)}}{m_n} - \frac{1}{n})^{2}= E_w (\frac{w_1^{(n)}}{m_n} - \frac{1}{n})^{2}=\frac{(1-\frac{1}{n})}{n m_n}  }.$$
We now employ Chebeshev's inequality  to bound $L_2(n)$ above as follows.
\begin{eqnarray*}
L_2(n)&\leq & \frac{m^{2}_{n}}{\ep^{\prime^{2}}  (1-\frac{1}{n})^2 } E_w\Big(\sum_{i=1}^{n} \big( \frac{w_i^{(n)}}{m_n} - \frac{1}{n})^{2}-\frac{(1-\frac{1}{n})}{m_n}    \Big)^{2}\\
&=& \frac{m^{2}_{n}}{\ep^{\prime^{2}}  (1-\frac{1}{n})^2 } \Big\{ E_w\Big(\sum_{i=1}^{n} \big( \frac{w_i^{(n)}}{m_n} - \frac{1}{n})^{2} \Big)^{2}- \frac{(1-\frac{1}{n})^2}{m^{2}_{n}}   \Big\}\\
&=& \frac{m^{2}_{n}}{\ep^{\prime^{2}}  (1-\frac{1}{n})^2 } \Big\{ n E_w \big( \frac{w^{(n)}_{1}}{m_n}-\frac{1}{n} \big)^{4}+ n(n-1) E_w \Big( \big(\frac{w^{(n)}_{1}}{m_n}-\frac{1}{n}\big)^2 \big(\frac{w^{(n)}_{2}}{m_n}-\frac{1}{n}\big)^2    \Big) - \frac{(1-\frac{1}{n})^2}{m^{2}_{n}}   \Big\}.
\end{eqnarray*}
 In view of the fact that  $w^{(n)}_{i}$, $1\leq i \leq n$ have multinomial distribution,          after computing  $E_{w} \big[(w^{(n)}_1)^{a} (w^{(n)}_2)^{b}\big]$, where $a,b$ are  two integers such that  $0\leq a,b\leq 2$, followed by some algebra,  we can bound the preceding term by
 \begin{eqnarray*}
&&\frac{m^{2}_{n}}{\ep^{\prime^{2}}  (1-\frac{1}{n})^2 } \Big\{  \frac{(1-\frac{1}{n})}{n^{3}m^{3}_{n}}+ \frac{(1-\frac{1}{n})^{4}}{m^{3}_{n}}+\frac{(m_n -1)(1-\frac{1}{n})^{2}}{n m^{3}_{n}}+ \frac{4 (n-1)}{n^{3} m_n}\\
  && +  \frac{1}{m^{2}_{n}}- \frac{1}{n m^{2}_n}+\frac{n-1}{n^{3}m^{2}_{n}} + \frac{4(n-1)}{n^{2} m^{3}_{n}}- \frac{(1-\frac{1}{n})^{2}}{m^{2}_n} \Big\}\\
 &&\sim \frac{1}{\ep^{\prime^{2}}}  \Big\{   \frac{4 m_n}{n^{2}}+  \frac{1}{n^3 m_n}+ \frac{1}{m_{n}}+  \frac{1}{n^{2}}+ \frac{4}{n m_n}      \Big\},
\end{eqnarray*}
where $a_n\sim b_n$ stands for the asymptotic equivalence of numerical sequences $a_n$ and $b_n$.
\par
Clearly, as $n,m_n\rightarrow \infty$,  the preceding relation approaches zero when $m_n=o(n^2)$. Now the  proof of Lemma \ref{lemma2} is complete. $\square$

\subsection*{Proof of Lemma \ref{lemma3}}
For $\ep_1,\ep_2,\ep_3 >0$ we have,
\begin{eqnarray*}
 &&P_w\big(   P_{X|w}\Big(\big|\frac{S_{m_n}^{*2}\Big/ m_n}{ \sigma^2\ \sumn\big( \frac{w_i^{(n)}}{m_n}-\frac{1}{n} \big)^2} -1 \big|>\ep_1\Big)>\ep_2 \big) \\
 &=&P_w\big(  \big\{ P_{X|w}\Big( \frac{\big| \frac{S_{m_n}^{*2}}{m_{n}}-\sigma^2 \sumn\big( \frac{w_i^{(n)}}{m_{n}}-\frac{1}{n} \big)^2\big|}{\sigma^2 \sumn\big(\frac{ w_i^{(n)}}{m_{n}}-\frac{1}{n} \big)^2}>\ep_1\Big)>\ep_2 \big)\\
&\leq &  P_w\big(  P_{X|w}\Big( \frac{\big| \frac{S_{m_n}^{*2}}{m_{n}}-\sigma^2 \sumn\big( \frac{w_i^{(n)}}{m_{n}}-\frac{1}{n} \big)^2\big|}{\sigma^2 \sumn\big(\frac{ w_i^{(n)}}{m_{n}}-\frac{1}{n} \big)^2}>\ep_1\Big)>\ep_2, \big| \frac{m_n}{(1-\frac{1}{n})} \sum_{i=1}^{n} \big( \frac{w_i^{(n)}}{m_n} - \frac{1}{n} \big)^2 -1 \big|\leq {\ep_3} \big)\\
 &+& P_w\big(   \sumn \big| \frac{m_n}{(1-\frac{1}{n})} \sum_{i=1}^{n} \big( \frac{w_i^{(n)}}{m_n} - \frac{1}{n} \big)^2 -1 \big| > {\ep_3}      \big)\\
&\leq &  P_w\big(   P_{X|w}\Big( \big| \frac{S_{m_n}^{*2}}{m_{n}}-\sigma^2 \sumn\big( \frac{w_i^{(n)}}{m_{n}}-\frac{1}{n} \big)^2\big|
 >\frac{\sigma^2 \ep_1(1-\ep_3)(1-\frac{1}{n})}{m_n}\Big)>\ep_2 \ \big)\\
&+& P_w\big(   \sumn \big| \frac{m_n}{(1-\frac{1}{n})} \sum_{i=1}^{n} \big( \frac{w_i^{(n)}}{m_n} - \frac{1}{n} \big)^2 -1 \big| > {\ep_3}     \big)\\
&=:& t_{1}(n)+t_{2}(n).
\end{eqnarray*}
We note  that along the lines of the proof of Lemma \ref{lemma2} it was already  shown that, when $m_n=o(n^2)$, as $n \to \infty$, then $t_{2}(n)\to 0$.
\par
To show that $t_{1}(n)\to 0$, as $n \to \infty$,  we proceed as follows.

\begin{eqnarray}
t_{1}(n)&\leq& P_w\big(  P_{X|\mathcal{G}_n}\Big( \big| \frac{S_{m_n}^{*2}}{m_{n}}- \frac{S^{2}_n}{m_{n}}   \big|>\frac{\sigma^2 \ep_1(1-\ep_3)(1-\frac{1}{n})}{m_n}\Big)>\frac{\ep_2}{3}  \big) \nonumber \\
&+&  P_w\big(  P_{X}\Big( \big|\frac{S^{2}_n}{m_{n}}-\frac{\sigma^2(1-\frac{1}{n})}{m_{n}} \big|>\frac{\sigma^2 \ep_1(1-\ep_3)(1-\frac{1}{n})}{m_n}\Big)>\frac{\ep_2}{3}  \big) \nonumber \\
&+& P_w\big(  I \Big( \big|   \frac{ \sum_{i=1}^{n}\big( \frac{w_i^{(n)}}{m_{n}}-\frac{1}{n} \big)^2- (1-\frac{1}{n})}{m_{n}} \big|>\frac{ \ep_1(1-\ep_3)(1-\frac{1}{n})}{m_n}\Big)>\frac{\ep_2}{3}  \big)\nonumber\\
&=:& t^{(1)}_{1}(n)+t^{(2)}_{1}(n)+t^{(3)}_{1}(n). \nonumber
%&=& \sum_{1\leq i\neq j\leq k} (\frac{w_i^{(k)}w_j^{(k)}}{m_k(m_k-1)}
%-\frac{1}{ k(k-1)})    (X_i-X_j)^2  \nonumber \\
%&+&  S^{2}_n-\sigma^{2} \nonumber\\
%&+& \sigma^{2}-\frac{\sigma^2}{m_n} \sumn\big( w_i^{(n)}-\frac{m_n}{n} \big)^2.\nonumber
\end{eqnarray}
$\vspace{-.3 cm}$
Now from  the $U$-statistic representation of the sample variance we have that
$$
S_{m_n}^{*2}- S^{2}_n   = \mathop{\sum\sum}_{1\leq i\neq j\leq n}
\big( \frac{w_i^{(n)}w_j^{(n)}}{m_n(m_n-1)}
- \frac{1}{ n(n-1)}\big) \frac{(X_i-X_j)^2}{2}.
$$
Therefore, $t^{(1)}_{1}(n)$ can be bounded above by
\begin{eqnarray}
&&P_w\big(  P_{X|w}\big( \big|\mathop{\sum\sum}_{1\leq i\neq j\leq n}
(\frac{w_i^{(n)}w_j^{(n)}}{m_n(m_n-1)}
- \frac{1}{ n(n-1)} ) (X_i-X_j)^2\big| > {\sigma^2 \ep_1(1-\ep_3)(1-\frac{1}{n})}\big)>\frac{\ep_2}{3}  \big)\nonumber \\[2ex]
&\leq & P_w\big(  E_{X|w}\Big( \big|\mathop{\sum\sum}_{1\leq i\neq j\leq n}
(\frac{w_i^{(n)}w_j^{(n)}}{m_n(m_n-1)}
- \frac{1}{ n(n-1)} ) (X_i-X_j)^2\big| \Big) > {\sigma^2 \ep_1(1-\ep_3)(1-\frac{1}{n})} \frac{\ep_2}{3}  \big)\nonumber \\[2ex] \nonumber
&\leq&  P_w\big(   \mathop{\sum\sum}_{1\leq i\neq j\leq n}
\big| \frac{w_i^{(n)}w_j^{(n)}}{m_n(m_n-1)}
- \frac{1}{ n(n-1)} \big| E_{X}(X_i-X_j)^2 > {\sigma^2 \ep_1(1-\ep_3)(1-\frac{1}{n})} \frac{\ep_2}{3}  \big)\nonumber \\[2ex]
&\leq&  P_w\big(   \mathop{\sum\sum}_{1\leq i\neq j\leq n}
\big| \frac{w_i^{(n)}w_j^{(n)}}{m_n(m_n-1)}
- \frac{1}{ n(n-1)} \big|  > { \ep_1(1-\ep_3)(1-\frac{1}{n})} \frac{\ep_2}{6}  \big).\nonumber
\end{eqnarray}
\par
For  ease of  notation we set $\ep_{n}:= { \ep_1(1-\ep_3)(1-\frac{1}{n})} \frac{\ep_2}{6}$. Using this,  the preceding term  can be bounded above by
$\vspace{-.1 cm}$
\begin{eqnarray*}
&& \ep^{-2}_{n}\Big\{ n(n-1) E_{w} \big(  \frac{w^{(n)}_1 w^{(n)}_2}{m_n (m_n-1)}-\frac{1}{n(n-1)} \big)^2 \\
&+& n(n-1)(n-2) E_{w}\Big( \big|\frac{w^{(n)}_1 w^{(n)}_2}{m_n (m_n-1)}-\frac{1}{n(n-1)}\big|   \big|\frac{w^{(n)}_1 w^{(n)}_3}{m_n (m_n-1)}-\frac{1}{n(n-1)}\big|   \Big)\\
&+& n(n-1)(n-2)(n-3) E_{w}\Big( \big|\frac{w^{(n)}_1 w^{(n)}_2}{m_n (m_n-1)}-\frac{1}{n(n-1)}\big|   \big|\frac{w^{(n)}_3 w^{(n)}_4}{m_n (m_n-1)}-\frac{1}{n(n-1)}\big|   \Big)  \Big\}
\\
&\leq& \ep^{-2}_{n} \Big\{ n(n-1) E_{w} \big(  \frac{w^{(n)}_1 w^{(n)}_2}{m_n (m_n-1)}-\frac{1}{n(n-1)} \big)^2\\
&+& n(n-1)(n-2) E_{w} \big(  \frac{w^{(n)}_1 w^{(n)}_2}{m_n (m_n-1)}-\frac{1}{n(n-1)} \big)^2 \\
&+& n(n-1)(n-2)(n-3)  E_{w} \big(  \frac{w^{(n)}_1 w^{(n)}_2}{m_n(m_n -1) } -\frac{1}{n(n-1)} \big)^2 \Big\}\\
&\sim& \ep^{-2}_{n} \Big\{ \frac{n^2}{n^{2} m^{2}_{n}  } + \frac{n^3}{n^2 m^{2}_{n} } + \frac{n^4}{n^2  m^{2}_{n}} \Big\}\longrightarrow 0. \label{S_n fixed proof}
\end{eqnarray*}
The preceding conclusion, which implies that $t^{(1)}_{1}(n)\to 0$, is true since, as $n \to \infty$, $\ep_n \to   \ep_1(1-\ep_3) \frac{\ep_2}{6}$ and $n=o(m_n)$ by assumption, as $n,m_n \to \infty$. Moreover, we note that the preceding convergence to $0$ also takes place when $n$, the number of the original observations, is fixed and $m_{n}:=m \to \infty$ (cf. Remark \ref{S_n fixed}).
\par
To show $t^{(2)}_{1}(n)\to 0$, as $n\to \infty$,  we note that
\begin{eqnarray}
t^{(2)}_{1} &\leq& P_w\big(  P_{X}\Big( \big| S^{2}_n-\sigma^2 \big|>\sigma^2 \ep_1(1-\ep_3)(1-\frac{1}{n}) \big\}  \Big)>\frac{\ep_2}{3} \big) \nonumber \\
&\leq& 3\ep^{-1}_{2} P_{X}\Big( \big| S^{2}_n-\sigma^2 \big|>\sigma^2 \ep_1(1-\ep_3)(1-\frac{1}{n}) \big\}  \Big) \longrightarrow 0. \nonumber
\end{eqnarray}
\par
To deal with $t^{(3)}_{1}(n)$, we observe that it can be bounded above by
$$ 3 \ep^{-1}_2 P_w \Big(  \big|   \frac{(1-\frac{1}{n})}{m_{n}}- \sum_{i=1}^{n}\big( \frac{w_i^{(n)}}{m_{n}}-\frac{1}{n} \big)^2 \big|>\frac{ \ep_1 (1-\frac{1}{n})}{m_n}  \Big).$$
Once again we note that during  the proof of Lemma \ref{lemma2} it was shown that when $m_n=o(n^2)$, as $n\to \infty$,    the preceding term approaches zero, i.e., $t^{(3)}_{1}(n) \to 0$. We now conclude that, as $n\to \infty$,  $t_{1}(n)\rightarrow 0$, and the latter also completes  the proof of   Lemma \ref{lemma3}. $\square$
\subsection*{Proof of Corollary \ref{cor2}}
 Recall that $m_n:= \sumn \zeta_i= n \frac{m_n}{n}= n \bar{\zeta}$. In view of Theorem \ref{thm1},  the proof of parts (a) and (b) of Corollary \ref{cor2} will  result from showing  that, as $n\to \infty$,
\begin{numcases}{ M_n=\frac{\max_{1\leq i \leq n}  (\frac{\zeta_i}{n\bar{\zeta}_{n}}- \frac{1}{n} )^2}{\sum_{i=1}^{n} (\frac{\zeta_i}{n\bar{\zeta}_{n}}- \frac{1}{n} )^2}=}
o(1)\ a.s.-P_{\zeta}\ when\ E_{\zeta} (\zeta_{1}^{4})<\infty & \\
o_{P_{\zeta}}(1)\ when\ E_{\zeta} (\zeta_{1}^{2})<\infty  . &
\end{numcases}
Since $\zeta_{i}$'s are positive random variables, we have
\begin{eqnarray*}
M_n&=&\frac{\max_{1\leq i \leq n}  (\frac{\zeta_i}{n\bar{\zeta}_{n}}- \frac{1}{n} )^2}{\sum_{i=1}^{n} (\frac{\zeta_i}{n\bar{\zeta}_{n}}- \frac{1}{n} )^2}\\
&=&\frac{\max_{1\leq i \leq n}  (\zeta_i-\bar{\zeta_n})^2}{\sum_{1\leq i \leq n}  (\zeta_i-\bar{\zeta_n})^2}.
\end{eqnarray*}
In view of  Kolmogorov's strong law of large numbers, when $E_{\zeta} (\zeta_{1}^{2})<\infty$, we have that, as $n\to \infty$,
$$
\sum_{i=1}^{n} (\zeta_i-\bar{\zeta_n})^2 \Big/ n \to \ var(\zeta_1)\ a.s.-P_{\zeta}.
$$
Also,
$$
\frac{\max_{1\leq i \leq n}  \big|\zeta_i-\bar{\zeta_n\big|}} {\sqrt{n}} \leq \frac{2  \max_{1\leq i \leq n}  \zeta_i} {\sqrt{n}}.
$$
Therefore, to prove parts (a) and (b) of  Corollary (\ref{cor2}), it suffices to, respectively,  show that, as $n \to \infty$,
\begin{numcases}{ \frac{\max_{1\leq i \leq n}  \zeta_i}{\sqrt{n}}=}
o(1)\ a.s.-P_{\zeta}\ when\ E_{\zeta} (\zeta_{1}^{4})<\infty \label{name1}& \\
o_{P_{\zeta}}(1)\ when\ E_{\zeta} (\zeta_{1}^{2})<\infty. \label{name2} &
\end{numcases}

\par
To establish (\ref{name1}), for $\ep>0$,  we write
\begin{eqnarray*}
\sum_{n=1}^{\infty} n \ P_{\zeta}\big(  \zeta_1 > \ep \sqrt{n } \big) &\leq& \sum_{n=1 }^{\infty} E_{\zeta} \big\{ \zeta_{1}^{2} I\big( |\zeta_{1}|> \ep \sqrt{n}  \big)  \big\}\\
&=& \sum_{k=1}^{\infty} \sum_{n=1}^{k} E_{\zeta} \big\{ \zeta_{1}^{2} I\big(\ep \sqrt{k} < \zeta_{1}\leq \ep \sqrt{k+1}  \big)  \big\}\\
&\leq&  \sum_{k=1}^{\infty} k \ E_{\zeta} \big\{ \zeta_{1}^{2} I\big(\ep \sqrt{k} < \zeta_{1}\leq \ep \sqrt{k+1}  \big)  \big\}\\
&\leq& \sum_{k=1}^{n}   E_{\zeta} \big\{ \zeta_{1}^{4} I\big( \ep \sqrt{k} < \zeta_{1}\leq \ep \sqrt{k+1}  \big)  \big\}\\
&=& \ep E_{\zeta}  \big(\zeta_{1}^{4}\big)<\infty.
\end{eqnarray*}
\par
In order to prove (\ref{name2}) for $\ep>0$,  we  continue as follows.
\begin{equation*}
n  P_{\zeta}\big(  \zeta_1  > \ep \sqrt{n } \big)
\leq \ep^{-2}  E_{\zeta} \big( \zeta^{2}_{1} I (\zeta_{1}> \ep \sqrt{n}  ) \big) \longrightarrow 0,\ as \ n\to \infty.
\end{equation*}
This also  completes the proof of (\ref{name2}) and that of Corollary \ref{cor2}. $\square$
\subsection*{ Proof of Theorem \ref{theorem2} and Theorem \ref{theorem3}}
We first prove  Theorem \ref{theorem3}.
$\vspace{-.7 cm}$
\subsection*{Proof of Theorem \ref{theorem3}}
In order to prove this theorem, first define

\begin{equation}\label{jadid}
T^*_{m_n}(\mu) := \frac{\ds\sumn\big(\frac{ w_i^{(n)}}{m_n} - \frac{1}{n}\big)(X_i-\mu)} {\sqrt{\ds
\ds \frac{(1-\frac{1}{n})}{n\ m_n}   \sumn (X_i-\mu)^2    }}.
\end{equation}
Recall that $(w_1^{(n)},\ldots,w_n^{(n)}) \substack{d\\=} \hbox{multinomial}(m_n,\frac{1}{n},\ldots,\frac{1}{n}
)$ for each $n\geq 1$.  Hence, by virtue of Corollary 4.1 of \cite{Morris}, conditioning on the data, $T^*_{m_n}(\mu)$ is a {\em properly normalized} linear function of $w_1^{(n)},\ldots,w_n^{(n)}$.  The term {\em properly normalized} is
used since
$$
\sumn \hbox{var}_{w|X}
\Big((\frac{w_i^{(n)}}{m_n} -  \frac{1}{n})\ (X_i-\mu)\Big)
=\sumn \frac{(X_i-\mu)^2}{nm_n}\big(1-\frac{1}{n}\big).
$$
The latter  normalizing sequence is that of  the CLT of Corollary 4.1 of \cite{Morris}.
\par
It will be first shown that the conditions under which Corollary 4.1 of Morris \cite{Morris} holds true are satisfied in probability $P_X$.  Then, by making use of the characterization of convergence in probability by almost sure convergence of subsequences,  it will be concluded that for each subsequence $\{n_\ell\}_\ell$ of $\{n\}^{\infty}_{n=1}$, there is a further subsequence $\{n_{\ell_s}\}_s$, along which, from Corollary 4.1 of Morris \cite{Morris}, $T^*_{m_n}(\mu)$,  conditionally on the sample, converges in distribution to standard normal a.s.$-P_X$. The latter means that, $\forall\, t\in \mathds{R}$
\begin{equation}\label{eq18}
 P_{w|X}(T^*_{m_n}(\mu)\leq t) \to P(Z\leq t) \ in \ probability-P_X,
\end{equation}
where $Z \ \substack{d\\=} \ N(0,1)$.

The  conditions of Corollary 4.1 of Morris \cite{Morris} are satisfied, for  one has
\begin{itemize}
\item[(a)] ${\frac{m_n}{n} \geq \ep > 0},$  assumed,
\item[(b)] ${\max_{1\leq i\leq n} \Big(\frac{1}{n}\Big) \to 0}$, as $n\to\infty$,
\item[(c)]$\ds{\frac{\max_{1\leq i\leq n} \ var_{w|X}
\big( ( \frac{w_i^{(n)}}{m_n} - \frac{1}{n}) (X_i-\mu)\big)}
{\sumn Var_{w|X}
\big(  (\frac{w_i^{(n)}}{m_n} -  \frac{1}{n} )(X_i-\mu)\big)}
=\frac{\max_{1\leq i\leq n}(X_i-\mu)^2}{\sumn (X_i-\mu)^2} } \to 0. $
\end{itemize}
The latter holds true in probability-$P_X$.
\par
Conclusion (c) is a  characterization of $X\in DAN$ (cf., e.g.,  \cite{GineGotzeMason}). In view of (a), (b) and (c),  one can conclude  that (\ref{eq18}) holds true.

Now observe that for $T^*_{m_n, S_n}$, as defined in (\ref{eqTSn}), we have
$$
T^*_{m_n,S_n} =
\frac{S_n}
{\sqrt{  \frac{(1-\frac{1}{n})}{n} \sumn(X_i-\mu)^2  }} \ T_{m_n}^*(\mu).
$$

Via Slutsky's theorem in probability$-P_X$, one will have, $\forall t\in \mathds{R}$, as $n,m_n \to \infty$ as in (\ref{Newly added})
\begin{equation}\label{eq19}
P_{w|X}(T^*_{m_n,S_n}\leq t )  \longrightarrow  P(Z\leq t) \hbox{ ~in }\ probability-P_X,
\end{equation}
if it is shown that, for $\ep_1,\ep_2>0$, as $n \to \infty$,
\begin{equation}\label{eqadded1primeprime}
P_X\Big(  P_{w|X}\big(\frac{1}{S^{2}_n}\big|  S_n^2-\frac{\sumn (X_i-\mu)^2}{n}\big| > \ep_1\big)>\ep_2\Big)\to 0.
\end{equation}
In order to prove the preceding result, for $X \in DAN$, without loss of generality we first assume that $\mu=0$ and write
\begin{eqnarray}
&& P_X\Big(  P_{w|X}\big(\frac{ \bar{X}^{2}_{n}}{S^{2}_n} > \ep_1\big)>\ep_2\Big)\nonumber\\[2ex]
&=& P_X\Big(I\big(\frac{\bar{X}_{n}^{2}}{S^{2}_n}> \ep_1\big)>\ep_2\Big)\nonumber\\[2ex]
&\leq & \ep_2^{-1} P_X\Big(\frac{\bar{X}_{n}^{2}}{S^{2}_n} > \ep_1
\Big) \to 0, \hbox{ as } n\to\infty.\nonumber
\end{eqnarray}
The preceding relation is due to the laws of large numbers when $E_{X} X^2<\infty$. When $E_{X} X^2 = \infty$, then we also make use of  Raikov's theorem  (cf., e.g., \cite{GineGotzeMason}).  Hence (\ref{eq19}) is valid, and  the proof of Theorem \ref{theorem3} is complete. $\square$
$\vspace{-.5 cm}$
\subsection*{Proof of Theorem \ref{theorem2}}
Due to Slutsky's theorem in probability$-P_X$, Theorem \ref{theorem2} will follow if one shows that, for $\ep>0$, as $n, m_n \to \infty$ as in (\ref{Newly added}), we have

\begin{equation}\label{S^2DAN}
P_{w|X}\big(\frac{1}{S_{n}^{2}}\big|S^{*^2}_{m_n} - S_n^2\big|>\ep_1\big) \to 0 \ in \ probability-P_X.
\end{equation}
\par
Using the $U$-statistic representation of the sample variance, for $S^{*^2}_{m_n}$ and $S_n^2$ we write
\begin{eqnarray*}
&&S^{2}_{n}= \frac{1}{n^{2}} \sum_{1\leq i \neq j \leq n} (X_i -X_j)^{2}= \frac{n-1}{n} .\frac{1}{2 n(n-1)} \sum_{1\leq i \neq j \leq n} (X_i -X_j)^{2}\\
&& S^{*^2}_{m_n}= \frac{1}{2 m_n (m_n-1)} \sum_{1\leq i \neq j \leq n} w_{i}^{(n)} w_{j}^{(n)}  (X_i -X_j)^{2}.
\end{eqnarray*}
To establish (\ref{S^2DAN}), we fist note that when $E_{X} X^2=+\infty$, as $n \rightarrow +\infty$, $S^{2}_{n}\big/ \ell^{2}(n)\rightarrow 1$ in $P_{X}$ and  for $\ep_1,\ep_2$ and $\ep_3>0$, on using the above   $U$-statistic representations, we have
$\vspace{-.3 cm}$
\begin{eqnarray}
\qquad&~&~~P_X\big\{
P_{w|X}\big(\frac{1}{\ell^{2}(n)}\big|S^{*^2}_{m_n} - S_n^2\big|>\ep_1\big)>\ep_2\big\}\nonumber \\[2ex]
&\leq& P_X\big\{
P_{w|X}\Big(\frac{\big|S^{*^2}_{m_n} - S_n^2\big|}{\ell^{2}(n)}>\ep_1,
\bigcap_{1\leq i\neq j\leq n} \big| \frac{w_{i}^{(n)} w_{j}^{(n)} }{m_n(m_n -1)} -\frac{1}{n^2}\big| \leq \frac{\ep_3}{n^2 \sqrt{\log n}}   \Big)>\frac{\ep_2}{2}\big\}
\nonumber\\[2ex]
&& \qquad + I\big\{ P_w\big(\bigcup_{1\leq i\neq j\leq n} \big| \frac{w_{i}^{(n)} w_{j}^{(n)} }{m_n(m_n -1)} -\frac{1}{n^2} \big| > \frac{\ep_3}{n^2 \sqrt{\log n}}  \big)>\frac{\ep_2}{2} \big\},\nonumber\\
&\leq& P_X\big\{ I\big(  \frac{\sum_{1\leq i\neq j \leq n} (X_{i}-X_{j})^{2}}{2 n^{2}\ell^{2}(n) \sqrt{\log n} }>\frac{\ep_{1}}{\ep_3} \big) >\frac{\ep_2}{2}    \big\} \nonumber \\
&&  \qquad + I\big\{ P_w\big(\bigcup_{1\leq i\neq j\leq n} \big| \frac{w_{i}^{(n)} w_{j}^{(n)} }{m_n(m_n -1)} -\frac{1}{n^2} \big| > \frac{\ep_3}{n^2 \sqrt{\log n} }  \big)>\frac{\ep_2}{2} \big\}\nonumber\\
&\leq&  P_X\big\{  \frac{ \sum_{1\leq i\neq j \leq n} (X_i - X_j)^{2}  }{ 2 n^2\ell^{2}(n) \sqrt{\log n}}> \frac{\ep_1 \ep_2}{2\ep_3}   \big\}\nonumber\\
&& \qquad + \frac{2}{\ep_2} P_w\big(\bigcup_{1\leq i\neq j\leq n} \big| \frac{w_{i}^{(n)} w_{j}^{(n)} }{m_n(m_n -1)} -\frac{1}{n^2} \big| > \frac{\ep_3}{n^2 \sqrt{\log n} }  \big)\nonumber\\
&\leq&  o(1)+ \frac{2}{\ep_2} n^2 \exp\big\{  - \frac{m_n (m_n-1)}{n^2  \log n } . \frac{\ep^{2}_3}{2(1+\frac{\ep_3}{ \sqrt{\log n} })} \big\}\label{eq10}\\
&=& o(1). \nonumber
\end{eqnarray}
The  relation (\ref{eq10}) is due to the fact that $X\in DAN$,  and  an application of Bernstien's inequality for $w_{i}^{(n)} w_{j}^{(n)}$,  viewed as $\sum_{1\leq s \leq m_n (m_n -1)} I\big( Y_s=1 \big),$ where, $Y_s$, $1\leq s \leq m_n (m_n-1)$, are i.i.d. random variables which are uniformly distributed on the set $\{1,\ldots,n^2 \}.$ And this completes the proof of (\ref{eqadded1}).
\par
In order to prove  (\ref{eqadded1prime}) we first note that, as $n\to \infty$,
\begin{equation*}
\frac{\max_{1\leq j \leq n} (X_j - \mu)^2}{\sumn (X_i-\mu)^2} \longrightarrow 0\ a.s.-P_X
\end{equation*}
once again from  Corollary 4.1 of Morris \cite{Morris}, on taking $m_n=n$ and as $n \to \infty$, for $T_{n}^{*}(\mu)$  as defined in (\ref{jadid})  we conclude that, for all $t\in \mathds{R}$,
\begin{equation*}
P_{w|X}(T_{n}^{*}(\mu)\leq t) \longrightarrow P(Z \leq t)\ a.s.-P_{X}.
\end{equation*}
Now, in view of  (\ref{eqadded1primeprime}) and Slutsky's theorem, the proof of (\ref{eqadded1prime}) follows if we show that for $\varepsilon_1, \varepsilon_2>0$, as  $n,m_n \to \infty$ such that $m_n/n \to \infty$,
\begin{equation*}
P_{X}\big\{ \limsup_{n\to \infty}P_{w|X} \big( \big|   S_{m_n}^{*^2}-S_{_{n}}^{2}  \big|>\varepsilon_1      \big)      >\varepsilon_2 \big\}=0
\end{equation*}
Observing  now that $P_{w|X} \big( \big|   S_{m_{n}}^{*^2}-S_{n}^{2}  \big|>\varepsilon_1      \big)$ asymptotically is bounded above by
\begin{equation*}
 \sum_{1\leq i<j \leq n} \varepsilon^{-1} E_{w}\big| \frac{w^{(n)}_{i} w^{(n)}_{i}}{m_n(m_n-1)} -\frac{1}{n^{2}}   \big|  \big| \frac{(X_i -X_j)^2 }{2}-\sigma^2 \big|,
\end{equation*}
where $\sigma^2=var_{X} (X)$. We note that
\begin{equation*}
 E_{w}\big| \frac{w^{(n)}_{i} w^{(n)}_{i}}{m_n(m_n-1)} -\frac{1}{n^{2}}   \big|
\leq \sqrt{var_{w}(\frac{w^{(n)}_{i} w^{(n)}_{i}}{m_n(m_n-1)})}
\sim \frac{1}{m_n n}.
\end{equation*}
Observing now that, as $n\rightarrow \infty$,  $n^{-2} \sum_{1\leq i<j \leq n} \big| \frac{(X_i -X_j)^2 }{2}-\sigma^2 \big|$ is  convergent  a.s.-$P_X$ and that $n/m_n \to 0$, completes the proof
 of (\ref{eqadded1prime}) and also that of Theorem \ref{theorem2}. $\square$

\subsection*{Proof of Theorem \ref{Efron confidence}}
   Observe that, as $n,m_n$ approach infinity,  the asymptotic equivalence  of $S_{m_{n}}^{*^{2}}(b)\big/m_n$, $S^{2}_{n}\sum_{i=1}^{n} \big(\frac{w^{(n)}_{i}(b)}{m_n}-\frac{1}{n}\big)^{2}$  and $ \sigma^2\ \sum_{i=1}^{n} \big(\frac{w^{(n)}_{i}(b)}{m_n}-\frac{1}{n}\big)^{2}$   with respect to the conditional probability  $P_{X|w}$, for each $1\leq b \leq B$, yields asymptotic  equivalence  for   $T^{*}_{m_n}(b)$, $T^{**}_{m_n}(b)$ and $T^{*}_{m_{n},\sigma}(b)$, where the latter  is defined by
\begin{equation}\nonumber
T^{*}_{m_{n},\sigma}(b):= \frac{\sum_{i=1}^{n} \big( \frac{w^{(n)}_{i}(b)}{m_{n}}-\frac{1}{n} \big)X_{i}  }{\sigma\ \sqrt{\sum_{i=1}^{n} \big(   \frac{w^{(n)}_{i}(b)}{m_{n}}-\frac{1}{n}    \big)^{2} }  }, \ 1\leq b \leq B.
\end{equation}
Therefore, we only give the proof of this theorem for $T^{*}_{m_{n},\sigma}$ and its associated bootstrapped quantile which is defined by
$$
C^{(B)}_{\sigma,\alpha}:=\inf\{t: \ \frac{1}{B}\sum_{ b=1}^{B} I(T^{*}_{m_{n},\sigma}(b)\leq t)\geq \alpha \}.
$$
In other words, we shall show that, as  $n,m_n,B\rightarrow \infty$, we have
$$
C^{(B)}_{\sigma,\alpha}\longrightarrow z_{\alpha} \ in \ probability - \bigotimes_{b=1}^{\infty}P_{X,w(b)}.
$$
To do so, we first note that in view of  the asymptotic normality of $T^{*}_{m_{n},\sigma}(b)$, for each $1\leq b\leq B$,  one can conclude the asymptotic conditional independence of $T^{*}_{m_{n},\sigma}(b) $ and $T^{*}_{m_{n},\sigma}(b^\prime)$ for each  $1\leq b\neq b^{\prime}\leq B$, from the fact that conditionally they are  asymptotically      uncorrelated. The latter is established in the following Lemma \ref{asymptotic independence}.

\begin{lemma}\label{asymptotic independence}
Assume the  conditions  of Theorem \ref{Efron confidence}. As $n,m_n\rightarrow \infty$, for each $1\leq b\neq b^{\prime}\leq B$, we have
$$ E\Big( T^{*}_{m_{n},\sigma}(b)\  T^{*}_{m_{n},\sigma}(b^\prime)      \big| \textbf{(}w^{(n)}_{1}(b),\ldots,w^{(n)}_{n}(b)\textbf{)}, \textbf{(}w^{(n)}_{1}(b^{\prime}),\ldots,w^{(n)}_{n}(b^{\prime})\textbf{)}  \Big)\rightarrow 0 \ a.s.\ P_{w}. $$
\end{lemma}
\subsection*{Proof of Lemma \ref{asymptotic independence}}

For  ease of  notation,  we let $E_{.|b}(.)$ and $E_{.|b,b^{\prime}}(.)$  be the respective  short hand notations  for the conditional expectations $E\Big(  ~ . ~    \big| \textbf{(}w^{(n)}_{1}(b),\ldots,w^{(n)}_{n}(b)\textbf{)} \Big)$  and  $E\Big(  ~ . ~    \big| \textbf{(}w^{(n)}_{1}(b),\ldots,w^{(n)}_{n}(b)\textbf{)}, \\ \textbf{(}w^{(n)}_{1}(b^{\prime}),\ldots,w^{(n)}_{n}(b^{\prime})\textbf{)}  \Big)$.  Similarly, we let $P_{.|b}(.)$ and $P_{.|b,b^{\prime}}(.)$ stand for the conditional probabilities $P\Big(  ~ . ~    \big| \textbf{(}w^{(n)}_{1}(b),\ldots,w^{(n)}_{n}(b)\textbf{)}  \Big)$  and  $P\Big(  ~ . ~    \big| \textbf{(}w^{(n)}_{1}(b),\ldots,w^{(n)}_{n}(b)\textbf{)}, \\ \textbf{(}w^{(n)}_{1}(b^{\prime}),\ldots,w^{(n)}_{n}(b^{\prime})\textbf{)}  \Big)$, respectively.

\par
Now observe that from the  independence  of the $X_{i}$'s,  we conclude that

$$
E_{X|b,b^{\prime}}\big(  T^{*}_{m_{n},\sigma}(b)\  T^{*}_{m_{n},\sigma}(b^\prime)     \big)
= \frac{\sum_{i=1}^{n}  \big( \frac{w^{(n)}_{i}(b)}{m_n}-\frac{1}{n} \big) \big( \frac{w^{(n)}_{i}(b^{\prime})}{m_n}-\frac{1}{n} \big)    } { \sqrt{\sum_{k=1}^{n}  \big( \frac{w^{(n)}_{k}(b)}{m_n}-\frac{1}{n} \big)^{2}}   \sqrt{\sum_{l=1}^{n}  \big( \frac{w^{(n)}_{l}(b)}{m_n}-\frac{1}{n} \big)^{2}}   }.
$$
By this, with  $\ep_1,\ep_2$ and $\ep_3>0$,  we can write
\begin{eqnarray*}
&&P\Big( \big| E_{X|b,b^{\prime}}\big(  T^{*}_{m_{n},\sigma}(b)\  T^{*}_{m_{n},\sigma}(b^\prime)     \big) \big|>\ep_1    \Big)\\
&\leq& P \Big(  \frac{m_n}{(1-\frac{1}{n}) }\big| \sum_{i=1}^{n}  \big( \frac{w^{(n)}_{i}(b)}{m_n}-\frac{1}{n} \big) \big( \frac{w^{(n)}_{i}(b^{\prime})}{m_n}-\frac{1}{n} \big)   \big|  >\ep_1(1-\ep_2)(1-\ep_3)    \Big)\\
&+& P \Big(  \big| \frac{m_n}{(1-\frac{1}{n})} \sumn \big( \frac{w^{(n)}_{i}(b)}{m_n}-\frac{1}{n}  \big)^{2}-1  \big|>\ep_2  \Big)\\
&+& P \Big(  \big| \frac{m_n}{(1-\frac{1}{n})} \sumn \big( \frac{w^{(n)}_{i}(b^{\prime})}{m_n}-\frac{1}{n}  \big)^{2}-1  \big|>\ep_3  \Big).
\end{eqnarray*}
The last two terms in the preceding relation have already been shown to approach zero as $\frac{m_n}{n^{2}}\to 0$. We now show that the first term approaches zero as well  in view of the following argument which relies on the facts that  $w^{(n)}_{i}$, $1\leq i \leq n$ are multinoialy distributed and that for each $1\leq i,j \leq n$, $w^{(n)}_{i}(b)$  and $w^{(n)}_{j}(b^{\prime})$ are i.i.d.'s  (in terms of $P_w$)  when $b\neq b^{\prime}$.
\par
In what will follow, for the ease of notation we put  $\ep_4:= \ep_1(1-\ep_2)(1-\ep_3) $.
\begin{eqnarray*}
&& P \Big(  \frac{m_n}{(1-\frac{1}{n}) }\big| \sum_{i=1}^{n}  \big( \frac{w^{(n)}_{i}(b)}{m_n}-\frac{1}{n} \big) \big( \frac{w^{(n)}_{i}(b^{\prime})}{m_n}-\frac{1}{n} \big)   \big|  >\ep_4    \Big)\\
&\leq& \ep^{-2}_4 \frac{m^{2}_{n}}{(1-\frac{1}{n})^2} \Big\{   n E^{2}_{b} \big( \frac{w^{(n)}_{1}(b)}{m_n}-\frac{1}{n} \big)^{2} + n(n-1) E^{2}_{b} \Big[ \big(\frac{w^{(n)}_{1}(b)}{m_n}-\frac{1}{n}\big) \big(\frac{w^{(n)}_{2}(b)}{m_n}-\frac{1}{n}\big)  \Big]      \Big\}\\
&=& \ep^{-2}_4 \frac{m^{2}_{n}}{(1-\frac{1}{n})^2} \Big\{   n\big( \frac{(1-\frac{1}{n})}{n m_n} \big)^{2}+
n(n-1) \big(  \frac{-1}{m_n n^{2}} \big)^{2}    \Big\}\\
&\leq& \ep^{-2}_4 \Big( \frac{1}{n}+\frac{1}{n^{2}(1-\frac{1}{n})^2  } \Big)\to 0.
\end{eqnarray*}

Now the proof of Lemma \ref{asymptotic independence} is complete. $\square$
\par
 We now continue the proof of Theorem \ref{Efron confidence} by showing that for any $\ep>0$, as $n,m_{n},B\rightarrow \infty$,
\begin{equation}\label{cs. Mas. method 1}
\bigotimes_{b=1}^{\infty}P_{X,w(b)}\big( C^{(B)}_{\sigma,\alpha}\leq z_{\alpha}-\ep  \big)\rightarrow 0,
\end{equation}

\begin{equation}\label{cs. Mas. method 2}
\bigotimes_{b=1}^{\infty}P_{X,w(b)}\big( C^{(B)}_{\sigma,\alpha}> z_{\alpha}+\ep  \big)\rightarrow 0.
\end{equation}
Observe that we have
\begin{equation}\label{cs. Mas. method 3}
\bigotimes_{b=1}^{\infty}P_{X,w(b)}\big( C^{(B)}_{\sigma,\alpha}\leq z_{\alpha}-\ep  \big)
\leq \bigotimes_{b=1}^{\infty}P_{X,w(b)}\big(    \frac{1}{B}\sum_{ b=1}^{B} I(T^{*}_{m_{n},\sigma}(b)\leq z_{\alpha}-\ep )\geq \alpha \big)
\end{equation}
and
\begin{equation}\label{cs. Mas. method 4}
\bigotimes_{b=1}^{\infty}P_{X,w(b)}\big( C^{(B)}_{\sigma,\alpha}> z_{\alpha}+\ep  \big)
\leq \bigotimes_{b=1}^{\infty}P_{X,w(b)}\big(    \frac{1}{B}\sum_{ b=1}^{B} I(T^{*}_{m_{n},\sigma}(b)\leq z_{\alpha}+\ep )< \alpha \big).
\end{equation}
In view of (\ref{cs. Mas. method 3}) and (\ref{cs. Mas. method 4}), the relations (\ref{cs. Mas. method 1}) and (\ref{cs. Mas. method 2}) will follow if for each $a\in \mathds{R}$, one shows that, as $n,m_n,B \to \infty$,
\begin{equation}\label{cs. Mas. method 5}
\bigotimes_{b=1}^{\infty}P_{X,w(b)}\big( \frac{1}{B}\sum_{ b=1}^{B} \big| I(T^{*}_{m_{n},\sigma}(b)\leq a )- \Phi(a)\big|>\ep \big)\rightarrow 0,
\end{equation}
where $\Phi(.)$ is the  standard normal distribution function.
\par
To establish (\ref{cs. Mas. method 5}),  we write
\begin{eqnarray*}
&& \bigotimes_{b=1}^{\infty}P_{X,w(b)}\big( \frac{1}{B}\sum_{ b=1}^{B} \big| I(T^{*}_{m_{n},\sigma}(b)\leq a )- \Phi(a) \big|>\ep \big)\\
&=& E \Big\{ P \big(  \frac{1}{B}\sum_{ b=1}^{B} \big| I(T^{*}_{m_{n},\sigma}(b)\leq a  )- \Phi(a)\big|>\ep \Big| {\bigotimes_{b=1}^{\infty}\mathfrak{F}_{w(b)}}      \big)     \Big\}\\
&\leq& E \Big\{ \frac{1}{B^2} \sum_{b=1}^{B} E_{X|b} \big( I(T^{*}_{m_{n},\sigma}(b)\leq a )- \Phi(a)\big)^{2} \Big\}\\
&+&  E \Big\{ \frac{1}{B^2} \sum_{1\leq b\neq b^{\prime}\leq B} E_{X|b,b^{\prime}}  \Big[ \big( I(T^{*}_{m_{n},\sigma}(b)\leq a )- \Phi(a)  \big)\ \big( I(T^{*}_{m_{n},\sigma}(b^{\prime})\leq a )- \Phi(a)  \big)         \Big]  \Big\}\\
&\leq& \frac{1}{B}+ E \Big\{  E_{X|1,2}  \Big[ \big( I(T^{*}_{m_{n},\sigma}(1)\leq a )- \Phi(a)  \big)\ \big( I(T^{*}_{m_{n},\sigma}(2)\leq a )- \Phi(a)  \big)         \Big]  \Big\}\\
&\longrightarrow& 0, \ as\ n,m_n,B\rightarrow \infty.
\end{eqnarray*}
The preceding relation is true since, in view of Lemma \ref{asymptotic independence}, for large enough $n,m_n$ we have that
\begin{eqnarray*}
&& E \Big\{ E_{X|1,2}  \Big[ \big( I(T^{*}_{m_{n},\sigma}(1)\leq a )- \Phi(a)  \big)\ \big( I(T^{*}_{m_{n},\sigma}(2)\leq a )- \Phi(a)  \big)         \Big]  \Big\}\\
& \approx& E\Big\{ E_{X|1} \big( I(T^{*}_{m_{n},\sigma}(1)\leq a )- \Phi(a)      \Big\}\ E\Big\{ E_{X|2} \big( I(T^{*}_{m_{n},\sigma}(2)\leq a )- \Phi(a)    \big)  \Big\}\\
&=& E \Big\{ P_{X|1} \big( T^{*}_{m_{n},\sigma}(1)\leq a\big) - \Phi(a)      \Big\}\
E \Big\{ P_{X|2} \big( T^{*}_{m_{n},\sigma}(2)\leq a \big)  - \Phi(a)      \Big\}\\
&&\longrightarrow 0, \ as \ n,m_n\rightarrow \infty.
\end{eqnarray*}
The preceding relation is due to part (a) of Corollary \ref{cor1}, with $\sigma^2$ replacing $S^{2}_{n}$ therein, and Lemma 1.2 in \cite{CsorgoRosalsky}. Now the proof of part (a) of Theorem \ref{Efron confidence} is complete.
\par
To prove parts (b) and (c), we first  conclude the  asymptotic in probability  equivalence  of $T^{**}_{m_{n}}$ and $T^{*}_{m_n,\mu}$, as $n,m_n\rightarrow \infty$, in terms of the conditional probability $P_{w|X}$ (cf. the proof of Theorem \ref{theorem2}). The same equivalence holds true  between $T^{*}_{m_n,S_n}$ and $T^{*}_{m_n,\mu}$ by virtue of Theorem \ref{theorem3} (cf. the proof of Theorem \ref{theorem3}).  Therefore, parts (b) and (c) will follow if we show that, as $n,m_n,B\rightarrow \infty$,
$$
C^{(B)}_{\mu,\alpha} \longrightarrow z_{\alpha}\ in \ probability - \bigotimes_{b=1}^{\infty}P_{X,w(b)},
$$
where $C^{(B)}_{\mu,\alpha}:= \inf\{t: \ \frac{1}{B}\sum_{ b=1}^{B} I(T^{*}_{m_{n},\mu}(b)\leq t)\geq \alpha \}$. To do so, similarly to what we did  in the proof of part (a), we shall show that for any $\ep>0$, as $n,m_n,B \to \infty$, we have
\begin{equation}\label{cs. Mas. method 10}
\bigotimes_{b=1}^{\infty}P_{X,w(b)} \big( C^{(B)}_{\mu,\alpha}\leq z_{\alpha}-\ep  \big)\rightarrow 0
\end{equation}
and
\begin{equation}\label{cs. Mas. method 11}
\bigotimes_{b=1}^{\infty}P_{X,w(b)} \big( C^{(B)}_{\mu,\alpha}> z_{\alpha}+\ep  \big)\rightarrow 0.
\end{equation}
Now observe that
\begin{equation}\label{cs. Mas. method 12}
\bigotimes_{b=1}^{\infty}P_{X,w(b)} \big( C^{(B)}_{\mu,\alpha}\leq z_{\alpha}-\ep  \big)
\leq \bigotimes_{b=1}^{\infty}P_{X,w(b)} \big(    \frac{1}{B}\sum_{ b=1}^{B} I(T^{*}_{m_{n},\sigma}(b)\leq z_{\alpha}-\ep )\geq \alpha \big)
\end{equation}
and
\begin{equation}\label{cs. Mas. method 13}
\bigotimes_{b=1}^{\infty}P_{X,w(b)} \big( C^{(B)}_{\mu,\alpha}> z_{\alpha}+\ep  \big)
\leq \bigotimes_{b=1}^{\infty}P_{X,w(b)}\big(    \frac{1}{B}\sum_{ b=1}^{B} I(T^{*}_{m_{n},\sigma}(b)\leq z_{\alpha}+\ep )< \alpha \big).
\end{equation}
In view of (\ref{cs. Mas. method 12}) and (\ref{cs. Mas. method 13}), the relations (\ref{cs. Mas. method 10}) and (\ref{cs. Mas. method 11}) will follow if for each $a\in \mathds{R}$, one shows that, as $n,m_n,B\to \infty$,
\begin{equation}\label{cs. Mas. method 14}
\bigotimes_{b=1}^{\infty}P_{X,w(b)}\big( \frac{1}{B}\sum_{ b=1}^{B} \big| I(T^{*}_{m_{n},\mu}(b)\leq a )- \Phi(a)\big|>\ep \big)\rightarrow 0.
\end{equation}
We establish the preceding relation in a similar way we established (\ref{cs. Mas. method 5})of part (a), on  noting that the proof here will be done via conditioning on the sample. Before sorting out the details, it is important to note that,  via conditioning on the sample,  $T^{*}_{m_{n},\mu}(b)$ and $T^{*}_{m_{n},\mu}(b^{\prime})$ are independent for each $1\leq b\neq b^{\prime} \leq B$. We have
\begin{eqnarray*}
&& \bigotimes_{b=1}^{\infty}P_{X,w(b)}\big( \frac{1}{B}\sum_{ b=1}^{B} \big| I(T^{*}_{m_{n},\mu}(b)\leq a )- \Phi(a) \big|>\ep \big)\\
&=& E \Big\{ P \big(  \frac{1}{B}\sum_{ b=1}^{B} \big| I(T^{*}_{m_{n},\mu}(b)\leq a )- \Phi(a)\big|>\ep  \Big| X   \big)     \Big\}\\
&\leq& E \Big\{ \frac{1}{B^2} \sum_{b=1}^{B} E \Big[  \Big( I(T^{*}_{m_{n},\mu}(b)\leq a )- \Phi(a)  \Big)^{2}      \Big| X \Big] \Big\}\\
&+&  E \Big\{ \frac{1}{B^2} \sum_{1\leq b\neq b^{\prime}\leq B} E  \Big[ \Big( I(T^{*}_{m_{n},\mu}(b)\leq a )- \Phi(a)  \Big)\ \Big( I(T^{*}_{m_{n},\mu}(b^{\prime})\leq a )- \Phi(a)  \Big)  \Big| X         \Big]  \Big\}\\
&\leq& \frac{1}{B}+ E \Big\{   \Big( P (T^{*}_{m_{n},\mu}(1)\leq a\big| X )- \Phi(a)  \Big)\ \Big( P (T^{*}_{m_{n},\mu}(2)\leq a \big| X )- \Phi(a)  \Big)           \Big\}\\
&&\longrightarrow 0, \ as\ n,m_n,B\rightarrow \infty.
\end{eqnarray*}
The preceding relation is true due to the fact that , as $n,m_n\rightarrow \infty$,
$$P(T^{*}_{m_{n},\mu}\leq a \big| X) \rightarrow \Phi(a)\ in \ probability - P_{X}$$
and Lemma 1.2 in S. Cs\"{o}rg\H{o} and Rosalsky \cite{CsorgoRosalsky}. Now the proof of (\ref{cs. Mas. method 14}) and,  consequently  that  of parts (b) and (c) are complete. Hence the proof of Theorem \ref{Efron confidence} is also complete.  $\square$
\subsection*{Proof of Theorem \ref{zeta confidence}}
Once again, in view of the fact that, as $n\rightarrow \infty$, $S^{2}_{n} \rightarrow \sigma^{2} a.s.-P_{X}$ we replace $T^{*}_{n}$ with $T^{*}_{n,\sigma}$,  which is defined by
$$
T^{*}_{n,\sigma}:= \frac{ \displaystyle \sum^n_{i=1} \big( \zeta_{i} - \bar{\zeta_n} \big) X_{i}}
{\sigma \sqrt{  \displaystyle\mathop{\sum}_{i=1}^{n}  (\zeta_{i}- \bar{\zeta}_n)^2} }\ a.s.-P_{\zeta}
$$
The proof of this theorem essentially consists of the same steps as those of part (a) of Theorem \ref{Efron confidence}. Hence, once again, the asymptotic normality of $T^{*}_{n,\sigma}(b)$, for each $1\leq b\leq B$,   conclude the asymptotic conditional independence of $T^{*}_{n,\sigma}(b) $ and $T^{*}_{n,\sigma}(b^\prime)$ for each  $1\leq b\neq b^{\prime}\leq B$, from the fact that conditionally they are   asymptotically   uncorrelated. The latter is established in the following Lemma \ref{asymptotic independence for zeta}.

\begin{lemma}\label{asymptotic independence for zeta}
Assume the conditions of Theorem \ref{zeta confidence}. As $n,m_n\rightarrow \infty$, for each $1\leq b\neq b^{\prime}\leq B$, we have that
$$ E\Big( T^{*}_{n,\sigma}(b)\  T^{*}_{n,\sigma}(b^\prime)      \big| \textbf{(}\zeta_{1}(b),\ldots,\zeta_{n}(b)\textbf{)}, \textbf{(}\zeta_{1}(b^{\prime}),\ldots,\zeta_{n}(b^{\prime})\textbf{)}  \Big)\rightarrow 0 \ in\ probability- P_{\zeta}. $$
\end{lemma}
To prove this lemma,  without loss of generality we assume that $E_{\zeta}(\zeta_1)=0$, and let $E_{.|b,b^{\prime}}$ be a short hand notation for  $E\Big( .\big| \textbf{(}\zeta_{1}(b),\ldots,\zeta_{n}(b)\textbf{)}, \textbf{(}\zeta_{1}(b^{\prime}),\ldots,\zeta_{n}(b^{\prime})\textbf{)}  \Big)$.
\par
Now,  similarly to the proof of  Lemma \ref{asymptotic independence for zeta}, we  note that
$$
E_{X|b,b^{\prime}}\big(  T^{*}_{n,\sigma}(b)\  T^{*}_{n,\sigma}(b^\prime)     \big)
= \frac{\sum_{i=1}^{n}  \big( \zeta_{i}(b)-\overline{\zeta}(b) \big) \big( \zeta_{i}(b)-\overline{\zeta}(b^{\prime}) \big)    } { \sqrt{ \sum_{k=1}^{n}  \big( \zeta_{k}(b)-\overline{\zeta}(b) \big)^{2} }  \sqrt{\sum_{l=1}^{n}  \big( \zeta_{l}(b)-\overline{\zeta}(b^{\prime}) \big)^{2}} }  .
$$
In view of the preceding statement, to complete the proof, with $\ep_1,\ep_2,\ep_3>0$, we  proceed as follows:
\vspace{.5 cm}
\begin{eqnarray*}
&& P\Big( \frac
{\big|  \sum_{i=1}^{n}  \big( \zeta_{i}(b)-\overline{\zeta}(b) \big) \big( \zeta_{i}(b^{\prime})-\overline{\zeta}(b^{\prime}) \big)\big|    } { \sqrt{ \sum_{k=1}^{n}  \big( \zeta_{k}(b)-\overline{\zeta}(b) \big)^{2} }  \sqrt{\sum_{l=1}^{n}  \big( \zeta_{l}(b)-\overline{\zeta}(b^{\prime}) \big)^{2}} }   >\ep_1\Big)\\
&\leq & P\Big( \frac{\big|  \sum_{i=1}^{n}  \big( \zeta_{i}(b)-\overline{\zeta}(b) \big) \big( \zeta_{i}(b^{\prime})-\overline{\zeta}(b^{\prime}) \big)\big|    } { n }   >\ep_1(1-\ep_2)(1-\ep_3) \Big)\\
&+& P\Big( \Big| \frac{\sum_{k=1}^{n}  \big( \zeta_{k}(b)-\overline{\zeta}(b) \big)^{2}}{n}  -1 \Big| >\ep_2\Big)\\
&+& P\Big( \Big| \frac{\sum_{k=1}^{n}  \big( \zeta_{k}(b^{\prime})-\overline{\zeta}(b^{\prime}) \big)^{2}}{n}  -1 \Big| >\ep_3\Big).
\end{eqnarray*}
Clearly, the last two relations  approach zero as $n\rightarrow \infty$. Hence, it only remains to show the asymptotic negligibility of the first term of the preceding three. To do so, we let $\ep_4:=\ep_1(1-\ep_2)(1-\ep_3)$ and apply Cheshev's inequality to arrive at
\begin{eqnarray*}
&& P\Big( \frac{\big|  \sum_{i=1}^{n}  \big( \zeta_{i}(b)-\overline{\zeta}(b) \big) \big( \zeta_{i}(b)-\overline{\zeta}(b^{\prime}) \big)\big|} {n}   >\ep_4 \Big)\\
&\leq& \ep^{-2}_4 n^{-2}  \big\{ n E^{2} \big( \zeta_{1}(b)-\overline{\zeta}(b) \big)   + n(n-1)   E^{2} \big[ \big( \zeta_{1}(b)-\overline{\zeta}(b) \big) \big( \zeta_{2}(b)-\overline{\zeta}(b) \big) \big]   \big\}\\
&\leq& \ep^{-2}_4 n^{-2} \big\{ n E^{2}(\zeta^{2}_1)+\frac{n(n-1)}{n^2} E^{2}(\zeta^{2}_1) \big\}\to 0, \ \textrm{as} \ n\rightarrow \infty.
\end{eqnarray*}
This completes the proof of Lemma \ref{asymptotic independence for zeta}. $\square$
\par
Due to similarity of the  rest of the proof of this theorem  and  that of (\ref{cs. Mas. method 5}) of part (a) in the proof of Theorem \ref{Efron confidence}, the details are omitted.  Now the proof of Theorem \ref{zeta confidence} is complete. $\square$


\begin{thebibliography}{9}

\bibitem{ArconesGine} \textsc{Arcones, M. A. and Gin\'{e}, E.} (1989). The Bootstrap of Mean with Arbitrary Bootstrap Sample Size. \textit{Annales de l'Institut Henri Poincar\'{e}} \textbf{25}, 1431-1452.


\bibitem{Arenal} \textsc{Arenal-Guti\'{e}rrez, E., Matr\'{a}n, C. and Cuesta-Albertos, J. A.} (1996). On the Unconditional Strong Law of Large Numbers for the Bootstrap Mean. \textit{Statistics and Probability Letters} \textbf{27}, 49-60.



\bibitem{Matran} \textsc{Arenal-Guti\'{e}rrez, E. and Matr\'{a}n, C.} (1996). A Zero-One Law Approach to the Central Limit Theorem for the Weighted Bootstrap Mean. \textit{Annals of Probability} \textbf{24}, 532-540.

\bibitem{Athreya} \textsc{Athrya, K. b.} (1987). Bootstrap of the mean in the infinite variance case
 \textit{The Annals of Statistics} \textbf{15},  724-731.

\bibitem{Billingsley} \textsc{Billingsley, P.} (1999). \textit{Probability and Measure}, 3nd ed.
Wiley, New York.

\bibitem{Csorgo and Nasari} \textsc{Cs\"{o}rg\H{o}, M. and Nasari M. M.} (2012). Asymptotics of Randomly Weighted $u$- and $v$-statistics: Application to Bootstrap. arXiv:1210.2757.

\bibitem{Csorgo} \textsc{Cs\"{o}rg\H{o}, S.} (1992). On the Law of Large Numbers for the Bootstrap Mean. \textit{Statistics and Probability Letters} \textbf{14}, 1-7.

%\bibitem{Csorgo and Mason} \textsc{Cs\"{o}rg\H{o}, S. and Mason, D. M.} ().

\bibitem{CsorgoRosalsky} \textsc{Cs\"{o}rg\H{o}, S. and Rosalsky A.} (2003). A Survey of Limit Laws for Bootstrapped sums. \textit{International Journal of Mathematics and Mathematical Sciences} \textbf{45}, 2835-2861.

\bibitem{CsorgoMason} \textsc{Cs\"{o}rg\H{o}, S. and Mason D. M.} (1989). Bootstrapping Empirical Functions. \textit{Annals of Statistics} \textbf{17}, 1447-1471.

\bibitem{DasGupta} \textsc{DasGupta, A.} (2008).  \textit{Asymptotic Theory of Statistics and Probability}. Springer Verlag, New York.


\bibitem{Efron} \textsc{Efron, B.} (1979). Bootstrap methods: another look at the jackknife. \textit{Annals of Statistics} \textbf{7}, 1-26.

\bibitem{EfronTibshirani}\textsc{Efron, B. and Tibshrani R.} (1993). \textit{An Introduction to the Bootstrap}. Chapman \& Hall, New York Londoan.

%\bibitem{Egorov} \textsc{Beals, R.} (2004). \textit{Analysis: An Introduction}.  Cambridge University Press

\bibitem{Gine} \textsc{Gin\'{e}, E.} (1996). \textit{Lectures on some aspects of the bootstrap}. Lectures on Probability Theory and Statistics.  \textsc{Giné, E.,  Grimmett, G. R. and  Saloff-Coste, L.} (1996). Ecole d’Et\'{e} de Probabilit\'{e}s de Saint-Flour XXVI-1996.

\bibitem{GineGotzeMason} \textsc{Gin\'{e}, E., G\"{o}tze, F. and Mason D. M.} (1997). When is the student $t$-statistic asymptotically Normal? \textit{Annals of Probability} \textbf{25}, 1514-1531.

\bibitem{Hajeck} \textsc{H\'{a}jek, J.} (1961).  Some extensions of the Wald-Wolfowitz-Noether theorem. \textit{Annals of Mathematical Statistics},\textbf{32}, 506-523.

\bibitem{Hallunconditional} \textsc{Hall, P.} (1986). On the Bootstrap and Confidence Intervals. \textit{Annals of Statistics} \textbf{14}, 1431-1452.


\bibitem{Hall} \textsc{Hall, P.} (1990). Asymptotics of the Bootstrap for Heavy-tailed distribution. \textit{Annals of Probability} \textbf{18}, 1342-1360.

\bibitem{Mason and Newton}  \textsc{Mason D. M.} and \textsc{Newton, M. A.}  (1992). A Rank Statistics Approach to the Consistency of a General Bootstrap.  \textit{Annals of Statistics} \textbf{ 20}, 1611-1624.

\bibitem{MasonShao} \textsc{Mason, D. M. and Shao Q.} (2001). Bootstrapping the Student $t$-statistic. \textit{Annals of Probability} \textbf{29}, 1435-1450.

\bibitem{Morris} \textsc{Morris, C.} (1975). Central Limit Theorems for Multinomial Sums. \textit{Annals of Statistics} \textbf{3}, 165-188.


\bibitem{Rubin} \textsc{Rubin, D. B.} (1981). The Baysian Bootstrap. \textit{Annals of Statistics} \textbf{9}, 130-134.

\bibitem{Weng} \textsc{Weng, C.} (1989). On a Second-Order Asymptotic Property of the Bayesian Bootstrap Mean. \textit{Annals of Statistics} \textbf{17}, 705-710.



\end{thebibliography}
\end{document}